\documentclass[aop,MSNbibl,seceqn,dvips]{arximspdf}
\usepackage{mathrsfs}

\doi{10.1214/12-AOP809} 
\volume{42}
\issue{3}
\pubyear{2014}
\firstpage{959}
\lastpage{993}

\makeatletter
\newtheorem{theorem}{Theorem}[section]
\newtheorem{thma}{Theorem}
\newtheorem{lemma}[theorem]{Lemma}
\newtheorem{proposition}[theorem]{Proposition}
\newproclaim{remark}{Remark}
\newproclaim{fact}[theorem]{Fact}
\newtheorem{question}[theorem]{Question}
\newtheorem{conjecture}[theorem]{Conjecture}

\def\r{\mathbb{R}}
\def\P{\mathbf{P}}
\def\E{\mathbf{E}}
\def\Q{\mathbf{Q}}
\def\z{\mathbb{Z}}

\def\ee{\mathrm{e}}
\def\d{ \mathrm{d}}
\makeatother

\begin{document}
\begin{frontmatter}

\title{The Seneta--Heyde scaling for the branching random~walk}
\runtitle{The Seneta--Heyde scaling for the branching random walk}

\begin{aug}
\author[A]{\fnms{Elie} \snm{Aidekon}\corref{}\ead[label=e1]{elie.aidekon@upmc.fr}}
\and
\author[A]{\fnms{Zhan} \snm{Shi}\ead[label=e2]{zhan.shi@upmc.fr}}
\runauthor{E. Aidekon and Z. Shi}
\affiliation{Universit\'e Pierre et Marie Curie Paris VI}
\address[A]{Laboratoire de Probabilit\'es UMR 7599\\
Universit\'e Paris VI\\
4 place Jussieu\\
F-75252 Paris Cedex 05\\
France\\
\printead{e1}\\
\phantom{E-mail:\ }\printead*{e2}} 
\end{aug}

\received{\smonth{2} \syear{2012}}
\revised{\smonth{9} \syear{2012}}

%
\begin{abstract}
We consider the boundary case (in the sense of Biggins and
Kyprianou~[\textit{Electron. J. Probab.} \textbf{10} (2005) 609--631] in
a one-dimensional super-critical branching random walk, and study the
additive martingale $(W_n)$.
We prove that, upon the system's survival, $n^{1/2} W_n$ converges in
probability, but not
almost surely, to a positive limit. The limit is identified as a
constant multiple of the almost
sure limit, discovered by Biggins and Kyprianou~[\textit{Adv. in Appl. Probab.} \textbf{36} (2004) 544--581], of the derivative martingale.
\end{abstract}

%
\begin{keyword}[class=AMS]
\kwd{60J80}
\kwd{60F05}
\end{keyword}
\begin{keyword}
\kwd{Branching random walk}
\kwd{Seneta--Heyde norming}
\kwd{additive martingale}
\kwd{derivative martingale}
\end{keyword}

\end{frontmatter}

\section{Introduction}
\label{s:intro}

We consider a discrete-time one-dimensional branching
random walk, whose distribution is governed by a point process $\Theta$
on the line. The system starts with an initial particle at the origin.
At time $1$, the particle dies, giving birth to a certain number of new
particles. These new particles form the particles at generation 1. They
are positioned according to the distribution of the point process~$\Theta$;
it is possible that several particles share a same position.
At time $2$, each of these particles dies, while giving birth to new
particles that are positioned (with respect to the birth place)
according to the distribution of $\Theta$. And the system goes on
according to the same mechanism. At each generation, we assume that
particles produce new particles independently of each other and of
everything up to that generation.

We denote by $(V(x), |x|=n)$ the positions of the particles at the
$n$th generation; so $(V(x), |x|=1)$ is distributed as the point
process $\Theta$. The family of random variables $(V(x))$ is usually
referred to as a branching random walk (Biggins~\cite{biggins10}).
Clearly, the number of particles in each generation forms a
Galton--Watson process. We always assume that this Galton--Watson
process is super-critical, so the system survives with positive probability.

Throughout the paper, we assume the following condition:
%
\begin{equation}
\E \biggl( \sum_{|x|=1} \ee^{-V(x)} \biggr) =1,\qquad
\E \biggl( \sum_{|x|=1} V(x) \ee^{-V(x)} \biggr)
=0. \label{cond-hab}
\end{equation}

 The branching random walk is then said to be in the boundary
case (Biggins and Kyprianou~\cite{biggins-kyprianou05}). Loosely
speaking, under some mild integrability conditions, an arbitrary
branching random walk can always be made to satisfy (\ref{cond-hab})
after a suitable linear transformation, as long as either the point
process $\Theta$ is not bounded from below, or if it is, $\E[\sum_{|x|=1} \mathbf{1}_{ \{ V(x) = \underline{m}\} }]<1$, where $\underline
{m}$ denotes the essential infimum of $\Theta$. More detailed
discussions on the nature of assumption (\ref{cond-hab}) can be found
in (the ArXiv version of) Jaffuel~\cite{jaffuel}.

It is immediately seen that under assumption $\E[\sum_{|x|=1} \ee
^{-V(x)}] =1$,
\[
W_n:= \sum_{|x|=n} \ee^{-V(x)},\qquad n
\ge0,
\]
is a martingale (with respect to its natural filtration). In
the literature, $(W_n)$ is referred to as the \textit{additive martingale}
associated with the branching random walk. Since $(W_n)$ is
nonnegative, it converges almost surely to a (finite) limit, which,
under assumption $\E[\sum_{|x|=1} V(x)\ee^{-V(x)}]=0$, turns out to be
$0$; see Biggins~\cite{biggins-mart-cvg}, Lyons~\cite{lyons}. In
particular, $\min_{|x|=n} V(x) \to\infty$ almost surely on the set of
nonextinction\footnote{In fact, according to Biggins~\cite
{biggins-lindley}, this holds as long as $\E[\sum_{|x|=1} \ee^{-V(x)}
] =1$.}.

Many of the discussions in this paper are trivial if the system dies
out. So let us introduce the conditional probability
\[
\P^* (\bullet):= \P( \bullet| \mbox{nonextinction}).
\]

Under (\ref{cond-hab}), since $W_n \to0$, $\P^*$-almost surely (and
$\P
$-almost surely), the martingale is not uniformly integrable. It is
natural to ask at which rate $W_n$ goes to 0; in the literature, this
concerns the Seneta--Heyde norming for $W_n$, referring to the pioneer
work on Galton--Watson processes by Seneta~\cite{seneta} and
Heyde~\cite
{heyde}. The study of the Seneta--Heyde norming for the branching
random walk in a general context [i.e., without assuming (\ref
{cond-hab})] goes back at least to Biggins and Kyprianou~\cite
{biggins-kyprianou96} and~\cite{biggins-kyprianou97}. It was an open
problem of Biggins and Kyprianou~\cite{biggins-kyprianou05} to study
the Seneta--Heyde norming under assumption (\ref{cond-hab}). This
problem was recently investigated in~\cite{yzpolymer}, under suitable
integrability conditions.

\renewcommand{\thethma}{\Alph{thma}}
\begin{thma}[(\cite{yzpolymer})]\label{thA}
Assume (\ref{cond-hab}). If there exists
$\delta>0$ such that\break
$\E[(\sum_{|x|=1}1)^{1+\delta}]<\infty$ and that
$\E[ \sum_{|x|=1} \ee^{-(1+\delta)V(x)}] + \E[ \sum_{|x|=1} \ee^{
\delta V(x)}]<\infty$,
then there exists a deterministic sequence
$(\lambda_n)$ of positive numbers with
$0< \liminf_{n\to\infty} {\lambda_n \over n^{1/2}} \le\limsup_{n\to
\infty} {\lambda_n \over n^{1/2}} <\infty$,
such that under $\P^*$,
%
\begin{equation}
\lambda_n W_n \to\mathscr{W}^* \qquad\mbox{in distribution},
\label{seneta-heyde}
\end{equation}
 where $\mathscr{W}^*>0$ is a positive random variable.
 \end{thma}

Let us make a brief description of the law of $\mathscr{W}^*$. Consider
the distributional equation for the nonnegative random variable $Z$
(excluding the trivial solution $Z=0$),
\[
\mathscr{L}_Z(t) = \E^* \biggl\{ \prod
_{|x|=1} \mathscr{L}_Z\bigl(t\ee^{-V(x)}
\bigr) \biggr\}\qquad \forall t\ge0,
\]
where $\mathscr{L}_Z(t):= \E^* (\ee^{-t Z})$ denotes the
Laplace transform of $Z$. Under assumption~(\ref{cond-hab}), it is
known (Liu~\cite{liu98}, Biggins and Kyprianou~\cite
{biggins-kyprianou05}) that the equation has a unique positive solution
(up to multiplication by a constant), denoted by $\mathscr{W}^*$. The
Laplace transform $\mathscr{L}_Z$ can be considered as a traveling wave
solution to a discrete F-KPP equation.

One may wonder whether $\lambda_n$ can be taken to be (a constant
multiple of) $n^{1/2}$ in~(\ref{seneta-heyde}). Our main result,
Theorem \ref{tmain} below, will tell us that the answer is yes.

The study of the additive martingale $W_n$ relies on analyzing another
fundamental martingale. Let us define
%
\begin{equation}
\label{defderiv} D_n:= \sum_{|x|=n} V(x)
\ee^{-V(x)},\qquad  n\ge0.
\end{equation}

Since $\E[\sum_{|x|=1} V(x)\ee^{-V(x)}]=0$,\vspace*{1pt} one can easily
check that $(D_n)$ is also a martingale, with $\E(D_n)=0$; it is
referred to in the literature as the \textit{derivative martingale}
associated with the branching random walk. Convergence of this new
martingale was studied by Biggins and Kyprianou~\cite
{biggins-kyprianou04}. In order to state their result, we introduce the
following integrability conditions:
%
\begin{eqnarray}
\label{cond-2} \E \biggl[\sum_{|x|=1}
V(x)^2\ee^{-V(x)} \biggr] &<&\infty,
\\
\E\bigl[X \log_+^2 X \bigr] &<&\infty,\qquad \E[\widetilde{X} \log_+
\widetilde{X}] <\infty, \label{cond-X}
\end{eqnarray}
where $\log_+ y:=\max\{ 0, \log y\}$ and $\log_+^2 y:=
(\log
_+ y)^2$ for any $y\ge0$, and
\[
X:=\sum_{|x|=1}\ee^{-V(x)},\qquad \widetilde{X}:=
\sum_{|x|=1} V(x)^+ \ee^{-V(x)},
\]
with $V(x)^+:= \max\{ V(x), 0\}$. Throughout the paper, we
assume (\ref{cond-hab}), (\ref{cond-2}) and~(\ref{cond-X}). We believe
that these assumptions are optimal for our results.


\begin{thma}[(Biggins and Kyprianou~\cite{biggins-kyprianou04})]\label{thB}
Assuming (\ref{cond-hab}), (\ref{cond-2}) and (\ref
{cond-X}), we have
%
\begin{equation}
D_n \to D_\infty,\qquad \mbox{$\P^*$-a.s.}, \label{Dn->}
\end{equation}
the limit $D_\infty>0$ having the distribution of $\mathscr
{W}^*$ in (\ref{seneta-heyde}). 
\end{thma}
(The positiveness of $D_\infty$ was proved in \cite
{biggins-kyprianou04} under slightly stronger assumptions. To see why
it is valid under current assumptions, we refer to Proposition A.3 of
\cite{elie}.)

It is worth mentioning that although $D_n$ is a signed martingale, its
limit $D_\infty$ is $\P^*$-almost surely positive.

Our main result is as follows.

\begin{theorem}
\label{tmain}
Assume (\ref{cond-hab}), (\ref{cond-2}) and (\ref{cond-X}). Under
$\P
^*$, we have
%
\begin{equation}
\lim_{n\to\infty} n^{1/2} W_n = \biggl(
{2\over\pi\sigma^2} \biggr)^{ 1/2} D_\infty\qquad\mbox{in
probability}, \label{main}
\end{equation}
where $D_\infty>0$ is the random variable in Theorem \ref{thB}, and
\[
\sigma^2:= \E \biggl[ \sum_{|x|=1}
V(x)^2 \ee^{-V(x)} \biggr] \in(0, \infty).
\]
\end{theorem}

The convergence in probability in Theorem \ref{tmain} is optimal: it
cannot be strengthened into almost sure convergence, as is shown in the
following theorem.

\begin{theorem}
\label{tmain2}
Assume (\ref{cond-hab}), (\ref{cond-2}) and (\ref{cond-X}). We have
\[
\limsup_{n\to\infty} n^{1/2} W_n = \infty,\qquad
\mbox{$\P^*$-a.s.}
\]
\end{theorem}

Let us say a few words about the proof of the theorems.

The first step in the proof of Theorem \ref{tmain} consists of
introducing a truncated version of the martingales $W_n$ and $D_n$,
denoted by $W_n^{(\alpha)}$ and $D_n^{(\alpha)}$, respectively, where
$\alpha\ge0$ is a positive parameter. The truncation argument can be
traced back to Harris~\cite{harris99}; we use it in the context of
conditional spines, following the formalism of Kyprianou~\cite
{kyprianou}. Roughly speaking (for a rigorous treatment of such
approximations, see Section \ref{sproof}), when $n\to\infty$,
\[
W_n^{(\alpha)} \approx W_n, \qquad D_n^{(\alpha)}
\approx c_0 D_n,
\]
where $c_0\in(0, \infty)$ is a constant depending only on
the law of $\Theta$. Moreover, $D_n^{(\alpha)}$ is a nonnegative
martingale, which allows us to define a new probability, $\Q^{(\alpha
)}$. The distribution of the branching random walk under $\Q^{(\alpha
)}$ is characterized by Biggins and Kyprianou~\cite
{biggins-kyprianou04} in the form of a spinal decomposition
(recalled\vadjust{\goodbreak}
as Fact \ref{pspine}). By means of a second moment argument, we prove
in Proposition \ref{pmoments} that under~$\Q^{(\alpha)}$,
\[
n^{1/2} {W_n^{(\alpha)} \over D_n^{(\alpha)}} \to \theta\qquad \mbox{in probability,}
\]
where $\theta\in(0, \infty)$ is a constant. Finally, in
Section \ref{sproof}, by taking $\alpha$ to be a large (but fixed)
constant, we come back to the probability $\P^*$, and prove that under
$\P^*$, $n^{1/2} {W_n \over D_n} \to c_0 \theta= ( {2\over\pi\sigma
^2} )^{1/2}$ in probability. Together with Theorem \ref{thB}, this yields
Theorem \ref{tmain}.

Theorem \ref{tmain2} is proved in Section \ref{snon-as} by studying
the minimal position in the branching random walk. The main ingredient
is a well-known spinal decomposition for the branching random walk
(Lyons~\cite{lyons}). As a by-product, we give a new proof, but under
assumptions we believe to be optimal, of the fact that $\liminf_{n\to
\infty} {1\over\log n} \min_{|x|=n} V(x) = {1\over2}$, $\P^*$-a.s.

The rest of the paper is as follows.

\begin{itemize}
\item In Section \ref{snew-pair}, we introduce a one-dimensional
random walk $(S_n)$ associated with the branching random walk, and
collect a few elementary properties of $(S_n)$.

\item Section \ref{sspine}: formalism of the truncation argument.

\item Section \ref{smoments}: proof of convergence in probability
of $n^{1/2} {W_n^{(\alpha)} \over D_n^{(\alpha)}}$ under $\Q
^{(\alpha)}$.

\item Section \ref{sproof}: proof of Theorem \ref{tmain}.

\item Section \ref{snon-as}: proof of Theorem \ref{tmain2}.

\item In Section \ref{squestions}, a few questions are raised for
further investigations.
\end{itemize}

Let us mention that our method allows us to prove the analogues of
Theorems~\ref{tmain} and \ref{tmain2} for the branching Brownian
motion. In fact, the main ingredients in our proof, namely the
truncation argument and spinal decompositions, are known in the case of
the branching Brownian motion. We prefer not to give any details on how
to make necessary modifications to obtain the analogues of Theorems~\ref
{tmain} and~\ref{tmain2} for the branching Brownian motion. These
modifications are more or less painless; moreover, the situation for
the branching Brownian motion is often neater than for the branching
random walk---for example, the analogue of the $h$-process whose
transition probabilities are given by (\ref{h-process}), is the
three-dimensional Bessel process, which is a well-studied stochastic
process in the literature. Instead, we close this paragraph with an
anecdotal remark: the pioneering work of \mbox{McKean}~\cite{mckean75} gives
an important motivation of the study of the branching Brownian motion
by connecting it to the Fisher--Kolmogorov--Petrovsky--Piscounov
(\mbox{F-KPP}) differential equation. Taking the almost sure limit of a
positive martingale (which is the analogue of the additive martingale
$W_n$), \mbox{McKean} claims that its Laplace transform, after a simple scale
change, gives a traveling wave solution to the \mbox{F-KPP} equation. There
turns out to be a flaw in the argument, pointed out by \mbox{McKean}~\cite
{mckean76}. Later on, Lalley and Sellke show in \cite{lalley-sellke}
that the almost sure limit studied in \cite{mckean75} actually is 0;
instead, they use another martingale (the analogue of the derivative
martingale $D_n$), and prove that its almost sure limit, which is
positive, has the Laplace transform as being a traveling wave solution.
Now that we know the two martingales (with the additive martingale
suitably normalized) have similar asymptotic behaviors in probability,
it becomes clear that the martingale limits studied by McKean~\cite
{mckean75} and by Lalley and Sellke~\cite{lalley-sellke} are a.s.
identical---if the additive martingale in McKean~\cite{mckean75} is
suitably normalized.

Throughout the paper, we use $a_n \sim b_n$ ($n\to\infty$) to denote
$\lim_{n\to\infty} {a_n\over b_n} =1$; the letter $c$ with subscript
denotes a finite and positive constant. We also adopt the notation
$\min_\varnothing:= \infty$, $\sum_\varnothing:= 0$ and $\prod_\varnothing:=1$.
For $x\in\r\cup\{ \infty\} \cup\{-\infty\}$, we write $x^+$ for
$\max\{ x, 0\}$.

\section{One-dimensional random walks}
\label{snew-pair}

This section collects some well-known material. We first
introduce a one-dimensional random walk associated with our branching
random walk, and then recall a few ingredients of fluctuation theory
for one-dimensional random walks.

\subsection{An associated one-dimensional random walk}
\label{subsSn}

Let $(V(x))$ be a branching random walk satisfying
(\ref
{cond-hab}) and (\ref{cond-2}). For any vertex $x$, we denote by $[\! [
\varnothing, x]\! ]$ the unique shortest path relating $x$ to the root
$\varnothing$, and $x_i$ (for $0\le i\le|x|$) the vertex on $[\! [
\varnothing, x]\! ]$ such that $|x_i|=i$. Thus, $x_0=\varnothing$ and
$x_{|x|}=x$. In words, $x_i$ (for $i<|x|$) is the ancestor of $x$ at
generation $i$. We also write $ ]\! ] \varnothing, x]\! ]:= [\! [
\varnothing, x]\! ] \setminus\{ \varnothing\}$.

The assumption $\E[\sum_{|x|=1} \ee^{-V(x)}]=1$ guarantees the
existence of an i.i.d. sequence of real-valued random variables $S_1$,
$S_2-S_1, S_3-S_2, \ldots,$ such that for any $n\ge1$ and any
measurable function $g\dvtx \r^n \to[0, \infty)$,
%
\begin{equation}
\E \biggl\{ \sum_{|x|=n} g\bigl(V(x_1),
\ldots, V(x_n)\bigr) \biggr\} = \E \bigl\{ \ee^{S_n}
g(S_1, \ldots, S_n) \bigr\}. \label{many-to-one}
\end{equation}

The law of $S_1$ is, according to (\ref{many-to-one}), given by
\[
\E\bigl[f(S_1)\bigr] = \E \biggl\{ \sum
_{|x|=1} \ee^{-V(x)} f\bigl(V(x)\bigr) \biggr\},
\]
for any measurable function $f\dvtx \r\to[0, \infty)$. Since
$\E
[\sum_{|x|=1} V(x)\ee^{-V(x)}]=0$, we have $\E(S_1)=0$. Let
%
\begin{equation}
\sigma^2:= \E\bigl[S_1^2\bigr] = \E \biggl
\{ \sum_{|x|=1} V(x)^2 \ee^{-V(x)}
\biggr\}. \label{sigma}
\end{equation}

Under (\ref{cond-hab}) and (\ref{cond-2}), we have
$0<\sigma
^2 <\infty$.

It is easy to prove (\ref{many-to-one}) by induction on $n$; see, for
example, Biggins and Kyprianou~\cite{biggins-kyprianou97}. The
presence
of the new random walk $(S_i)$ is explained via a
change-of-probabilities technique as in Lyons, Pemantle and Peres~\cite
{lyons-pemantle-peres},\vadjust{\goodbreak} and Lyons~\cite{lyons}; see Fact \ref{fLyons}
for more details. In the literature, the change-of-probabilities
technique is used by many authors in various forms (see \cite
{lyons-pemantle-peres} for a detailed account), the idea going back at
least to Kahane and Peyri\`ere~\cite{kahane-peyriere}.

\subsection{Elementary properties of one-dimensional random walks}
\label{subsappendix}

Let $S_1,\break  S_2-S_1, S_3-S_2, \ldots$ be an
i.i.d. sequence of real-valued random variables with $\E(S_1) =0$ and $\sigma
^2:= \E[S_1^2] \in(0, \infty)$. Let $\tau^+:= \inf\{ k\ge1\dvtx S_k
\ge
0\}$, which is well defined almost surely (because $\E(S_1)=0$). Let
%
\begin{equation}
R(u):= \E \Biggl\{ \sum_{j=0}^{\tau^+ -1} {
\mathbf{1}}_{ \{ S_j \ge-u\} } \Biggr\},\qquad u\ge0, \label{h0}
\end{equation}
which, according to the duality lemma, is the renewal
function associated with the entrance of $(-\infty, 0)$ by the walk
$(S_n)$. More precisely, the function $R$ can be expressed as
%
\begin{equation}
\label{defh0} R(u) = \sum_{k=0}^\infty\P\bigl\{
|H_k| \le u\bigr\},\qquad u\ge0,
\end{equation}
where $H_0<H_1<H_2<\cdots$ are the strict descending ladder
heights of $(S_n)$; that is, $H_k:= S_{\tau^{-}_k}$, with $\tau^{-}_0:= 0$ and $\tau^{-}_k:= \inf\{ i> \tau^{-}_{k-1}\dvtx S_i < \min_{0\le
j\le\tau^{-}_{k-1}} S_j\}$, $k\ge1$.

Throughout the paper, we regularly use the following identity:
%
\begin{equation}
R(u) = \E \bigl\{ R(S_1 +u) \mathbf{1}_{ \{ S_1 \ge-u\} } \bigr\}\qquad \forall
u\ge0. \label{harmonic}
\end{equation}

Conditions $\E[S_1^2]<\infty$ and $\E(S_1)=0$ ensure that $\E(
|H_1| )
<\infty$; see, for example, \cite{feller}, Theorem XVIII.5.1. The
renewal theorem states that the limit
%
\begin{equation}
c_0:= \lim_{u\to\infty} {R(u) \over u}
\label{c0}
\end{equation}
exists and lies in $(0, \infty)$. As a consequence, there
exist constants $c_2 \ge c_1>0$ such that
%
\begin{equation}
c_1(1+u)\le R(u) \le c_2 (1+u),\qquad u\ge0. \label{h0>}
\end{equation}


The function $R(\cdot)$ describes the persistency of $(S_i)$. In fact,
if we write
\[
\underline{S}_n:= \min_{1\le i\le n} S_i,\qquad  n
\ge1,
\]
then there exists a constant $0<\theta<\infty$ such that
%
\begin{equation}
\P\{ \underline{S}_n \ge0 \} \sim {\theta\over n^{1/2}},\qquad n\to
\infty. \label{kozlov}
\end{equation}

More generally, for any $u\ge0$,
%
\begin{equation}
\P\{ \underline{S}_n \ge-u \} \sim {\theta R(u) \over n^{1/2}},\qquad  n\to
\infty. \label{kozlov2}
\end{equation}
See Kozlov~\cite{kozlov}, formula (12).\vadjust{\goodbreak}

We will need a uniform version of (\ref{kozlov2}) for $u$ depending on
$n$. Let $(b_n)$ be a sequence of positive numbers such that $\lim_{n\to\infty} {b_n\over n^{1/2}}=0$.\vspace*{1pt} Then (see \cite{elie-bruno}) for
any bounded continuous function $f\dvtx [0, \infty) \to\r$, we have, as
$n\to\infty$,
%
\begin{equation}
\qquad\E \biggl\{ f \biggl( {S_n +u \over(n\sigma^2)^{1/2}} \biggr) \mathbf{1}_{ \{ \underline{S}_n \ge-u \} }
\biggr\} = {\theta R(u) \over n^{1/2}} \biggl( \int_0^\infty
f(t) t\ee^{-t^2/2} \,\d t + o(1) \biggr), \label{elie-bruno}
\end{equation}
uniformly in $u\in[0, b_n]$. In particular,
%
\begin{equation}
\P\{ \underline{S}_n \ge-u \} \sim {\theta R(u) \over n^{1/2}}, \qquad n\to
\infty, \label{kozlov-uniform}
\end{equation}
uniformly in $u\in[0, b_n]$.

\begin{lemma}
\label{lrelation-constants}
Let $c_0$ and $\theta$ be the constants in
(\ref{c0}) and (\ref{kozlov}), respectively.
Then
%
\begin{equation}
\theta c_0 = \biggl( {2\over\pi\sigma^2} \biggr)^{ 1/2}.
\label{relation-constants}
\end{equation}
\end{lemma}

\begin{pf}
We recall from (\ref{defh0}) that $R(u)$ is the
mean number of strict descending ladder heights within $[-u, 0]$. By
the renewal theorem (see Feller~\cite{feller}, Section~XI.1), we have
$c_0 = {1\over\E( |H_1| )}$. On the other hand (Feller~\cite{feller},
Theorem XII.7.4),
\[
\sum_{n\ge1} s^n \P\{
\underline{S}_n \ge0 \} = \exp \biggl(\sum
_{n\ge1} {s^n\over n} \P\{ S_n \ge0 \}
\biggr).
\]
Since $\E(S_1)=0$ and $\E(S_1^2) <\infty$, it follows from
Theorem XVIII.5.1 of Feller~\cite{feller} that $c:= \sum_{n\ge1}
{1\over n} [\P\{ S_n \ge0 \} - {1\over2}]$ is well defined,
satisfying\break  $\E( |H_1| ) = {\sigma\over2^{1/2}} \ee^c$. Accordingly,
\[
\sum_{n\ge1} s^n \P\{
\underline{S}_n \ge0 \} \sim {\ee^{c}\over(1-s)^{1/2}},\qquad  s\uparrow1.
\]
By a Tauberian theorem (Feller~\cite{feller}, Theorem
XIII.5.5), this yields that
\[
\P\{ \underline{S}_n \ge0 \} \sim {\ee^c \over(\pi n)^{1/2}},\qquad  n\to
\infty.
\]

Comparing with (\ref{kozlov}), we get $\theta= {\ee^c
\over
\pi^{1/2}} = ({2 \over\pi\sigma^2})^{1/2} \E(|H_1|) = ({2 \over
\pi
\sigma^2})^{1/2} {1\over c_0}$, proving Lemma \ref
{lrelation-constants}.
\end{pf}

\begin{lemma}
\label{la1}
There exists
$c_3>0$ such that for
$u>0$, $a\ge0$, $b\ge0$ and $n\ge1$,
\[
\P \{ \underline{S}_n \ge-a, b-a\le S_n \le b-a+u \} \le
c_3 {(u+1)(a+1)(b+u+1)\over n^{3/2}}.
\]
\end{lemma}

\begin{pf}The inequality is proved in \cite{ezsimple} for
a certain value of $u$, say $1$; hence, the inequality holds for $u<1$.
The case $u>1$ boils down to the case $u\le1$ by splitting the
interval $[b-a, b-a+u]$ into intervals of lengths $\le1$, the number
of these intervals being less than $(u+1)$.
\end{pf}

\begin{lemma}
\label{la2}
There exists $c_4>0$ such that for $a\ge0$,
\[
\sup_{n\ge1} \E \bigl[ |S_n| \mathbf{1}_{\{ \underline{S}_n \ge-a\} } \bigr]
\le c_4 (a+1).
\]
\end{lemma}

\begin{pf} We need to check that for some $c_5>0$, $\E[
S_n \mathbf{1}_{\{ \underline{S}_n \ge-a\} } ] \le c_5 (a+1)$, $\forall
a\ge0$, $\forall n\ge1$.

Let $\tau_a^-:= \inf\{ i\ge1\dvtx S_i <-a\}$. Then $\{ \underline{S}_n
\ge-a\} = \{ \tau_a^- >n\}$; thus\break  $\E[ S_n \mathbf{1}_{\{ \underline
{S}_n \ge-a\} } ] = - \E[ S_n \mathbf{1}_{\{ \tau_a^- \le n\} } ]$,
which, by the optional sampling theorem, equals $\E[ (-S_{\tau_a^-})
\mathbf{1}_{\{ \tau_a^- \le n\} } ]$. Therefore, $\sup_{n\ge1} \E[ S_n
\mathbf{1}_{\{ \underline{S}_n \ge-a\} } ] = \E[ (-S_{\tau_a^-})]$.

It remains to check that $\E[ (-S_{\tau_a^-}) -a] \le c_6 (a+1)$ for
some $c_6>0$ and all $a\ge0$, under the assumption $\E( S_1^2)<\infty
$.\footnote{Assuming $\E(|S_1|^3)<\infty$, even more is true
(Mogulskii~\cite{mogulskii73}): we have $\sup_{a\ge0} \E[ (-S_{\tau
_a^-}) -a] <\infty$.} By a known trick (Lai~\cite{lai}) using the
sequence of strict descending ladder heights $0=:H_0<H_1<H_2<\cdots$, it
boils down to proving that $\E[ (-H_{\tau_H(-a)}) -a] \le c_7 (a+1)$
for some $c_7>0$ and all $a\ge0$, where $H_1, H_2-H_1, H_3-H_2,
\ldots,$ are i.i.d. \textit{negative} random variables with $\E
(|H_1|) <
\infty$, and $\tau_H(-a):= \inf\{ i\ge1\dvtx H_i < -a\}$. This, however,
is a special case of (2.6) of Borovkov and Foss~\cite
{borovkov-foss}.
\end{pf}

\begin{lemma}
\label{la3}
Let $0<\lambda<1$. There exists
$c_8>0$ such that for
$a,b\ge0$, $0\le u\le v$ and $n\ge1$,
%
\begin{eqnarray}\label{ezsimple-bis}
&&\P \Bigl\{ \underline{S}_{\lfloor\lambda n\rfloor} \ge-a, \min_{i\in[\lambda n, n]\cap\z}
S_i \ge b-a, S_n \in[ b-a+u,b-a+v] \Bigr\}
\nonumber
\\[-8pt]
\\[-8pt]
\nonumber
&&\qquad\le c_8 {(v+1)(v-u+1)(a+1)\over n^{3/2}}.
\end{eqnarray}
\end{lemma}

\begin{pf} We treat $\lambda n$ as an integer. Let $\P
_{\mbox{\scriptsize(\ref{ezsimple-bis})}}$ denote the probability
expression on the left-hand side of (\ref{ezsimple-bis}). Applying
the\vadjust{\goodbreak}
Markov property at time $\lambda n$, we see that $\P_{\mbox
{\scriptsize
(\ref{ezsimple-bis})}} = \E[ \mathbf{1}_{\{ \underline{S}_{\lambda n}
\ge
-a, S_{\lambda n} \ge b-a \} } f(S_{\lambda n}) ]$, where $f(r):= \P
\{ \underline{S}_{n-\lambda n} \ge b-a-r, S_{n-\lambda n} \in
[b-a-r+u,b-a-r+v]\}$ (for $r\ge b-a$). By Lemma \ref{la1}, $f(r) \le
c_3 {(v+1)(v-u+1)(a+r-b+1)\over n^{3/2}}$ (for $r\ge b-a$). Therefore,
\[
\P_{\mbox{\scriptsize(\ref{ezsimple-bis})}} \le {c_3(v+1)(v-u+1)\over n^{3/2}} \E\bigl[ (S_{\lambda n} +
a-b+1) \mathbf{1}_{\{
\underline{S}_{\lambda n} \ge-a, S_{\lambda n} \ge b-a \} } \bigr].
\]
The expectation $\E[\cdots]$ on the right-hand side being
bounded by $\E[ |S_{\lambda n}|\times  \mathbf{1}_{\{ \underline{S}_{\lambda n}
\ge-a\} } ] + a +1$, it suffices to apply Lemma \ref{la2}.
\end{pf}

\begin{lemma}
\label{la5}
There exists a constant $C>0$ such that for any
sequence $(b_n)$ of nonnegative numbers with
$\limsup_{n\to\infty} {b_n \over n^{1/2}}<\infty$,
and any $0<\lambda<1$,
we have
\[
\liminf_{n\to\infty} \inf_{b\in[0, b_n]} n^{3/2}
\P \Bigl\{ \underline{S}_{\lfloor\lambda n\rfloor} \ge0, \min_{\lfloor\lambda n\rfloor<j\le n}
S_j \ge b, b \le S_n \le b + C \Bigr\} > 0.
\]
\end{lemma}

\begin{pf} The lemma is proved in \cite{ezsimple} in the
special cases $\lambda= {1\over2}$ and $b=b_n$; the same proof is
valid for the general case $0<\lambda<1$ and uniformly in $b\in[0,
b_n]$.
\end{pf}

\begin{lemma}
\label{la6}
There exists a constant $c_9>0$ such that for any $y\ge0$ and $z\ge0$,
\[
\sum_{k\ge0} \P \{ S_k \le y-z,
\underline{S}_k \ge-z \} \le c_9 (1+y) \bigl(1+ \min\{y,
z\}\bigr).
\]
\end{lemma}

\begin{pf} See Lemma B.2(i) of \cite{elie}.
\end{pf}

\section{Truncated processes, change of probabilities}
\label{sspine}

In the study of the martingales $W_n$ and $D_n$, it
turns out to be more convenient to work with a truncated version of the
branching random walk. The truncating argument, originating from
Harris~\cite{harris99}, was formalized for the branching Brownian
motion in the context of the spine conditioned to stay positive by
Kyprianou~\cite{kyprianou}, and was later put into the branching random
walk setting by Biggins and Kyprianou~\cite{biggins-kyprianou04}. It
can be adapted in other situations, for example, in the study of
fragmentation processes (Bertoin and Rouault~\cite{bertoin-rouault},
Berestycki, Harris and Kyprianou~\cite{berestycki-harris-kyprianou}).

Let $(V(x))$ be a branching random walk. For any vertex $x$, we define
\[
\underline{V}(x):= \min_{y\in]\! ] \varnothing, x]\! ]} V(y).
\]

Let $\alpha\ge0$, and let $R(\cdot)$ be as in (\ref{h0}). Let
\[
R_\alpha(u):= R(u+\alpha),\qquad u\ge-\alpha.
\]

 Having in mind the additive martingale $(W_n)$ and the
derivative martingale $(D_n)$, let us introduce a new pair of processes
\begin{eqnarray*}
W_n^{(\alpha)} &:=& \sum_{|x|=n}
\ee^{-V(x)} \mathbf{1}_{\{ \underline{V}(x) \ge-\alpha\} },
\\
D_n^{(\alpha)} &:=&\sum_{|x|=n}
R_\alpha\bigl(V(x)\bigr) \ee^{-V(x)} \mathbf{1}_{\{ \underline{V}(x) \ge-\alpha\} }.
\end{eqnarray*}

Recall from (\ref{c0}) that $\lim_{u\to\infty} {R(u)\over u}
= c_0$. Under (\ref{cond-hab}), we have\break  $\inf_{|x|=n} V(x) \to\infty$,
$\P^*$-a.s. So, it is intuively clear that if $\alpha$ is ``sufficently
large,'' then $W_n^{(\alpha)}$ should behave like $W_n$, and
$D_n^{(\alpha)}$ like $c_0 D_n$. This can easily be made rigorous, and
will be done in Section \ref{sproof}.

In Section \ref{smoments}, we are going to prove that for any $\alpha
\ge0$, as $n\to\infty$, $n^{1/2} {W_n^{(\alpha)}\over D_n^{(\alpha)}}
\to\theta$ in probability [$\theta$ being the constant in (\ref
{kozlov})], under a new probability called $\Q^{(\alpha)}$. To define
this new probability $\Q^{(\alpha)}$, we first need a simple property
of $D_n^{(\alpha)}$. For any $n$, let $\mathscr{F}_n$ denote the
sigma-algebra generated by the branching random walk in the first $n$
generations.

The following result is known, and its analogue for the branching
Brownian motion is in \cite{kyprianou}.

\begin{fact}[(Biggins and Kyprianou~\cite{biggins-kyprianou04})]
Assume (\ref{cond-hab}). For any
$\alpha\ge0$, $(D_n^{(\alpha)}, n\ge0)$
is a nonnegative martingale with respect to
$(\mathscr{F}_n)$, such that
$\E(D_n^{(\alpha)}) = R_\alpha(0)$,
$\forall n$.
\end{fact}

Since $(D_n^{(\alpha)})$ is a nonnegative martingale with $\E
(D_n^{(\alpha)}) = R_\alpha(0)$, there exists a probability measure
$\Q
^{(\alpha)}$ such that for any $n$,
\[
\Q^{(\alpha)} |_{\mathscr{F}_n}:= {D_n^{(\alpha)} \over R_\alpha(0)} \bullet
\P\Big|_{\mathscr{F}_n}.
\]

We observe that $\Q^{(\alpha)}(\mbox{nonextinction})=1$,
and that $\Q^{(\alpha)} (D_n^{(\alpha)} >0) =1$ for any~$n$.

(Strictly speaking, to make our presentation mathematically rigorous,
we need to work on the canonical space of branching random walks ($={}$space of marked trees) and use the rigorous language of Neveu~\cite
{neveu86} to describe the probabilities $\P$ and $\Q^{(\alpha)}$, as
well as the forthcoming spine $(w^{(\alpha)}_n, n\ge0)$. We continue
using the informal language, and referring the interested reader to
Lyons~\cite{lyons} or Lyons and Peres~\cite{lyons-peres}, for a
rigorous treatment. We mention that in the next paragraph, while
introducing the spine $(w^{(\alpha)}_n)$, we should, strictly speaking,
enlarge the probability space and work on a product space.)

Recall that the positions of the particles in the first generation,
$(V(x),\break  |x|=1)$, are distributed under $\P$ as the point process
$\Theta$. Fix $\alpha\ge0$. For any real number $u\ge-\alpha$, let
$\widehat
{\Theta}_u^{(\alpha)}$ denote a point process whose distribution is the
law of $(u+V(x), |x|=1)$ under $\Q^{(u+\alpha)}$.

We now consider the distribution of the branching random walk under $\Q
^{(\alpha)}$. The system starts with one particle, denoted by
$w^{(\alpha)}_0$, at position $V(w^{(\alpha)}_0)=0$. At each step $n$
(for $n\ge0$), particles of generation $n$ die, while giving birth to
point processes independently of each other: the particle $w^{(\alpha
)}_n$ generates a point process distributed as $\widehat\Theta
^{(\alpha
)}_{V(w^{(\alpha)}_n)}$, whereas any particle $x$, with $|x|=n$ and
$x\not= w^{(\alpha)}_n$, generates a point process distributed as $V(x)
+ \Theta$. The particle $w^{(\alpha)}_{n+1}$ is chosen among the
children $y$ of $w^{(\alpha)}_n$ with probability proportional to
$R_\alpha(V(y)) \ee^{-V(y)} \mathbf{1}_{\{ \underline{V}(y) \ge-\alpha
\}
}$. The line of descent $w^{(\alpha)}:=(w^{(\alpha)}_n,n\ge0)$ is
referred to as the \textit{spine}. We denote by $\mathcal{B}^{(\alpha)}$
the family of the positions of this system.\footnote{The spine process
$w^{(\alpha)}$ is, of course, part of the new system. Since working in
a product space and dealing with projections and marginal laws would
make the notation complicated, we feel free, by a slight abuse of
notation, to identify $\mathcal{B}^{(\alpha)}$ with $(\mathcal
{B}^{(\alpha)}, w^{(\alpha)})$.}

\begin{fact}[(Biggins and Kyprianou~\cite{biggins-kyprianou04})]
\label{pspine}
Assume (\ref{cond-hab}). Let $\alpha\ge0$.
\begin{longlist}[(iii)]
\item[(i)] The branching random walk under $\Q^{(\alpha)}$,
has the distribution of $\mathcal{B}^{(\alpha)}$.

\item[(ii)] For any $n$ and any vertex $x$ with $|x|=n$, we have
%
\begin{equation}
\Q^{(\alpha)} \bigl\{ w^{(\alpha)}_n =x |
\mathscr{F}_n \bigr\} = {R_\alpha(V(x)) \ee^{-V(x)}
\mathbf{1}_{\{ \underline{V}(x) \ge
-\alpha\} } \over D_n^{(\alpha)}}. \label{wn}
\end{equation}

\item[(iii)] The spine process
$(V(w^{(\alpha)}_n), n\ge0)$ under
$\Q^{(\alpha)}$, is
distributed as the centered random
walk $(S_n, n\ge0)$ under $\P$
conditioned to stay in
$[-\alpha, \infty)$.
\end{longlist}
\end{fact}

Since $D_n^{(\alpha)}>0$, $\Q^{(\alpha)}$-a.s., identity (\ref{wn})
makes sense $\Q^{(\alpha)}$-almost surely. In Fact~\ref{pspine}(iii),
the centered random walk $(S_n)$ (under $\P$) conditioned to stay in
$[-\alpha, \infty)$ is in the sense of Doob's $h$-transform: it is a
Markov chain with transition probabilities given by
%
\begin{equation}
p^{(\alpha)}(u, \d v):= \mathbf{1}_{ \{ v\ge-\alpha\} } {R_\alpha(v) \over R_\alpha(u)}
p(u, \d v),\qquad u\ge-\alpha, \label{h-process}
\end{equation}
where $p(u, \d v):= \P(S_1 + u \in\d v)$ is the transition
probability of $(S_n)$. Fact \ref{pspine}(iii) tells that for any
$n\ge1$ and any measurable function $g\dvtx \r^n \to[0, \infty)$,
%
\begin{eqnarray}\label{spine}
&&\E_{\Q^{(\alpha)}} \bigl[g\bigl(V\bigl(w^{(\alpha)}_i\bigr), 0
\le i\le n\bigr)\bigr]
\nonumber
\\[-8pt]
\\[-8pt]
\nonumber
&&\qquad = {1\over R_\alpha(0)} \E \bigl[ g(S_i, 0
\le i\le n) R_\alpha(S_n) \mathbf{1}_{ \{ \underline{S}_n
\ge-\alpha\} } \bigr].
\end{eqnarray}


The spine decomposition will allow us, in the next section, to handle
the first two moments of ${W_n^{(\alpha)}\over D_n^{(\alpha)}}$ under
$\Q^{(\alpha)}$.

\section{\texorpdfstring{Convergence in probability of ${W_n^{(\alpha)}\over D_n^{(\alpha)}}$ under $\Q^{(\alpha)}$}
{Convergence in probability of W n (alpha)/D n (alpha) under Q (alpha)}}
\label{smoments}
The aim of this section is to prove that ${W_n^{(\alpha
)}\over D_n^{(\alpha)}}$ converges in probability (under $\Q^{(\alpha
)}$). We do this by estimating $\E_{\Q^{(\alpha)}} ({W_n^{(\alpha
)}\over D_n^{(\alpha)}})$ and $\E_{\Q^{(\alpha)}} [ ({W_n^{(\alpha
)}\over D_n^{(\alpha)}})^2]$, using Fact \ref{pspine} and its
consequence (\ref{spine}). Recall that $a_n \sim b_n$ ($n\to\infty$)
means $\lim_{n\to\infty} {a_n\over b_n} =1$.

\begin{proposition}
\label{pmoments}
Assume (\ref{cond-hab}), (\ref{cond-2}) and (\ref{cond-X}). Let
$\alpha
\ge0$.
We have
%
\begin{eqnarray}
\E_{\Q^{(\alpha)}} \biggl( {W_n^{(\alpha)}\over D_n^{(\alpha)}} \biggr) &\sim&
{\theta\over n^{1/2}}, \label{moment1}
\\
\E_{\Q^{(\alpha)}} \biggl[ \biggl( {W_n^{(\alpha)}\over D_n^{(\alpha)}}
\biggr)^{ 2} \biggr] &\sim& {\theta^2\over n},\qquad  n\to\infty,
\label{moment2}
\end{eqnarray}
where $\theta\in(0, \infty)$ is
the constant in (\ref{kozlov}). As
a consequence, under $\Q^{(\alpha)}$,
\[
\lim_{n\to\infty} n^{1/2} {W_n^{(\alpha)}\over D_n^{(\alpha)}} = \theta\qquad
\mbox{in probability}.
\]
\end{proposition}

The last part (convergence in probability) of the proposition is
obviously a consequence of (\ref{moment1})--(\ref{moment2}) and
Chebyshev's inequality.

The rest of the section is devoted to the proof of (\ref{moment1}) and
(\ref{moment2}). The first step is to represent ${W_n^{(\alpha)}\over
D_n^{(\alpha)}}$ as a conditional expectation. Recall that $\mathscr
{F}_n$ is the sigma-algebra generated by the first $n$ generations of
the branching random walk.

\begin{lemma}
\label{lE[WD|Fn]}
Assume (\ref{cond-hab}). Let
$\alpha\ge0$. We have, for any $n$,
\[
{W_n^{(\alpha)}\over D_n^{(\alpha)}} = \E_{\Q^{(\alpha)}} \biggl( {1\over R_\alpha(V(w^{(\alpha)}_n))} \Big|
\mathscr{F}_n \biggr),
\]
where $w^{(\alpha)}_n$ is, as before, the element
of the spine in the $n$th
generation.
\end{lemma}

\begin{pf}We have $\E_{\Q^{(\alpha)}} ( {1\over
R_\alpha
(V(w^{(\alpha)}_n))} | \mathscr{F}_n) = \sum_{|x|=n} {\Q^{(\alpha
)} \{
w^{(\alpha)}_n =x | \mathscr{F}_n \} \over R_\alpha(V(x))}$, which,
according to (\ref{wn}), equals $\sum_{|x|=n} {\ee^{-V(x)} \over
D_n^{(\alpha)}} \mathbf{1}_{\{\underline{V}(x)\ge-\alpha\} }=
{W_n^{(\alpha)}\over D_n^{(\alpha)}}$.
\end{pf}

We are now able to prove the first part of Proposition \ref{pmoments},
concerning $\E_{\Q^{(\alpha)}} ({W_n^{(\alpha)}\over
D_n^{(\alpha)}})$.

\begin{pf*}{Proof of Proposition \ref{pmoments}: Equation
(\ref
{moment1})} By Lemma \ref{lE[WD|Fn]}, $\E_{\Q^{(\alpha)}}\times  (
{W_n^{(\alpha)}\over D_n^{(\alpha)}}) = \E_{\Q^{(\alpha)}}
({1\over
R_\alpha(V(w^{(\alpha)}_n))})$, which,\vspace*{1pt} by applying (\ref{spine}) to
$g(u_0, u_1, \ldots, u_n):= {1\over R_\alpha(u_n)}$, equals ${\P\{
\underline{S}_n \ge-\alpha\}\over R_\alpha(0)}$. By (\ref{kozlov2}),\vspace*{2pt}
$\P\{ \underline{S}_n \ge-\alpha\} \sim{\theta R_\alpha(0) \over
n^{1/2}}$ (as $n\to\infty$), from which~(\ref{moment1}) follows
immediately.
\end{pf*}

It remains to prove (\ref{moment2}), which is done in several steps.
The first step gives the correct order of magnitude of $\E_{\Q
^{(\alpha
)}} [({W_n^{(\alpha)}\over D_n^{(\alpha)}} )^2]$:

\begin{lemma}
\label{lE[WD^2]}
Assume (\ref{cond-hab}) and (\ref{cond-2}). Let
$\alpha\ge0$. We have
\[
\E_{\Q^{(\alpha)}} \biggl[ \biggl( {W_n^{(\alpha)}\over D_n^{(\alpha)}}
\biggr)^{ 2} \biggr] = O \biggl( {1\over n} \biggr),\qquad n\to
\infty.
\]
\end{lemma}

\begin{pf} By Lemma \ref{lE[WD|Fn]} and Jensen's
inequality,
\begin{eqnarray*}
&&\E_{\Q^{(\alpha)}} \biggl[ \biggl( {W_n^{(\alpha)}\over D_n^{(\alpha
)}}
\biggr)^{ 2} \biggr] \\
&&\qquad\le \E_{\Q^{(\alpha)}} \biggl(
{1\over[R_\alpha(V(w^{(\alpha)}_n))]^2} \biggr).
\end{eqnarray*}
 The expression on the right-hand side is, by (\ref{spine}),
\begin{eqnarray*}
 &=&{1\over R_\alpha(0)} \E \biggl( {\mathbf{1}_{\{ \underline{S}_n \ge
-\alpha
\} }\over R_\alpha(S_n)} \biggr)\\
& =&
{1\over R_\alpha(0)} \E \biggl( {\mathbf{
1}_{\{ \underline{S}_n \ge-\alpha\} }\over R(S_n+\alpha)} \biggr).
\end{eqnarray*}

 Recall from (\ref{h0>}) that $R(u) \ge c_1(1+u)$, $\forall
u\ge0$. Therefore,
\begin{eqnarray*}
&&R_\alpha(0) c_1 \times \E_{\Q^{(\alpha)}} \biggl[ \biggl(
{W_n^{(\alpha)}\over D_n^{(\alpha)}} \biggr)^{ 2} \biggr]
\\
&&\quad\le \E \biggl( {\mathbf{1}_{\{ \underline{S}_n \ge-\alpha\} }
\over S_n+\alpha+1} \biggr)
\\
&&\quad\le \sum_{i=0}^{\lfloor n^{1/2}\rfloor-1} \E \biggl(
{\mathbf{1}_{\{ -\alpha+ i \le S_n
< -\alpha+ i+1,
\underline{S}_n \ge-\alpha\} }
\over S_n+\alpha+1} \biggr) + \E \biggl( {\mathbf{1}_{\{ S_n\ge-\alpha+ \lfloor n^{1/2}\rfloor,
\underline{S}_n \ge-\alpha\} }
\over S_n+\alpha+1} \biggr),
\end{eqnarray*}
 which, by Lemma \ref{la1}, is
\begin{eqnarray*}
&\le& \sum_{i=0}^{\lfloor n^{1/2}\rfloor-1}
{1\over i +1} c_3 {(\alpha+1) (i+1)\over n^{3/2}} +
{\P\{ \underline{S}_n \ge-\alpha\} \over
\lfloor n^{1/2}\rfloor}
\\
&=& {\lfloor n^{1/2}\rfloor c_3 (\alpha+1)
\over n^{3/2}} + {\P\{ \underline{S}_n \ge-\alpha\} \over
\lfloor n^{1/2}\rfloor}.
\end{eqnarray*}

 By (\ref{kozlov2}), $\P\{ \underline{S}_n \ge-\alpha\} =
O({1\over n^{1/2}})$, $n\to\infty$. The lemma follows.
\end{pf}

Lemma \ref{lE[WD^2]} tells us that $\operatorname{Var}_{\Q^{(\alpha)}}
( {W_n^{(\alpha)}\over D_n^{(\alpha)}} ) = O ( {1\over n} )$, whereas
our goal is to replace $O( {1\over n} )$ by $o( {1\over n} )$. We need
to do some more work.

Let $E_n$ be an event such that $\Q^{(\alpha)}(E_n) \to1$, $n\to
\infty
$. Let
\[
\xi_{n,E_n^c}:= \E_{\Q^{(\alpha)}} \biggl( {\mathbf{1}_{E_n^c} \over
R_\alpha(V(w^{(\alpha)}_n))} \Big|
\mathscr{F}_n \biggr).
\]

 Since ${W_n^{(\alpha)}\over D_n^{(\alpha)}} = \E_{\Q
^{(\alpha
)}} ( {1\over R_\alpha(V(w^{(\alpha)}_n))} | \mathscr{F}_n ) = \xi
_{n,E_n^c} + \E_{\Q^{(\alpha)}} ( {\mathbf{1}_{E_n} \over R_\alpha
(V(w^{(\alpha)}_n))} | \mathscr{F}_n )$, we have
\[
\E_{\Q^{(\alpha)}} \biggl[ \biggl( {W_n^{(\alpha)}\over D_n^{(\alpha)}}
\biggr)^{ 2} \biggr] = \E_{\Q^{(\alpha)}} \biggl[ {W_n^{(\alpha)}\over D_n^{(\alpha)}}
\xi_{n,E_n^c} \biggr] + \E_{\Q^{(\alpha)}} \biggl[ {W_n^{(\alpha)}\over D_n^{(\alpha)}}
{\mathbf{1}_{E_n}
\over R_\alpha(V(w^{(\alpha)}_n))} \biggr].
\]

 By the Cauchy--Schwarz inequality, we have
\begin{eqnarray*}
\E_{\Q^{(\alpha)}} \biggl[ {W_n^{(\alpha)}\over D_n^{(\alpha)}} \xi_{n,E_n^c} \biggr]
&\le& \biggl\{ \E_{\Q^{(\alpha)}} \biggl[ \biggl( {W_n^{(\alpha)}\over D_n^{(\alpha)}}
\biggr)^{ 2} \biggr] \biggr\}^{ 1/2} \bigl\{ \E_{\Q^{(\alpha)}}
\bigl(\xi_{n,E_n^c}^2\bigr)\bigr\}^{1/2}
\\
&=& O \biggl( {1\over n^{1/2}} \biggr) \bigl\{ \E_{\Q^{(\alpha)}}\bigl(
\xi_{n,E_n^c}^2\bigr)\bigr\}^{1/2},
\end{eqnarray*}
the last identity being a consequence of Lemma \ref
{lE[WD^2]}. So (\ref{moment2}) will be a straightforward
consequence of the following lemmas.

\begin{lemma}
\label{lxi}
Assume (\ref{cond-hab}) and
(\ref{cond-2}). Let
$\alpha\ge0$. For any sequence of
events $(E_n)$ such that
$\Q^{(\alpha)}(E_n) \to1$, we have
\[
\E_{\Q^{(\alpha)}} \bigl(\xi_{n,E_n^c}^2\bigr) = o \biggl(
{1\over n} \biggr),\qquad n\to\infty.
\]
\end{lemma}

\begin{lemma}
\label{lindependance}
Assume (\ref{cond-hab}),
(\ref{cond-2}) and (\ref{cond-X}). Let
$\alpha\ge0$. There exists a sequence
of events $(E_n)$ such that
$\Q^{(\alpha)}(E_n) \to1$, and that
\[
\E_{\Q^{(\alpha)}} \biggl[ {W_n^{(\alpha)}\over D_n^{(\alpha)}} {\mathbf{1}_{E_n} \over R_\alpha(V(w^{(\alpha)}_n))}
\biggr] \le{\theta^2 \over n} + o \biggl( {1\over n} \biggr),\qquad
n\to\infty.
\]
\end{lemma}

\begin{pf*}{Proof of Lemma \ref{lxi}} By Jensen's inequality,
\[
\E
_{\Q^{(\alpha)}}\bigl(\xi_{n,E_n^c}^2\bigr) \le\E_{\Q^{(\alpha)}}\biggl ( {\mathbf{
1}_{E_n^c} \over[R_\alpha(V(w^{(\alpha)}_n))]^2} \biggr).
\]
Consequently, for
any $\varepsilon>0$,
\begin{eqnarray*}
&&\E_{\Q^{(\alpha)}} \bigl(\xi_{n,E_n^c}^2\bigr)\\
&&\quad\le
\E_{\Q^{(\alpha)}} \biggl( {\mathbf{1}_{E_n^c} \over[R_\alpha(V(w^{(\alpha)}_n))]^2} \mathbf{1}_{\{ V(w^{(\alpha)}_n) \ge\varepsilon n^{1/2}\} } \biggr)
+ \E_{\Q^{(\alpha)}} \biggl( {\mathbf{1}_{\{ V(w^{(\alpha)}_n) < \varepsilon n^{1/2}\} }
\over[R_\alpha(V(w^{(\alpha)}_n))]^2} \biggr)
\\
&&\quad= \E_{\Q^{(\alpha)}} \biggl( {\mathbf{1}_{E_n^c} \over[R_\alpha(V(w^{(\alpha)}_n))]^2} \mathbf{1}_{\{ V(w^{(\alpha)}_n) \ge\varepsilon n^{1/2}\} }
\biggr) + \E \biggl( {\mathbf{1}_{\{ S_n < \varepsilon n^{1/2}\} }
\over R_\alpha(S_n) R_\alpha(0)} \mathbf{1}_{\{ \underline{S}_n\ge-\alpha\} } \biggr),
\end{eqnarray*}
the last identity being a consequence of (\ref{spine}).
Recall from (\ref{h0>}) that $R_\alpha(u) = R(u+\alpha) \ge c_1
(1+u+\alpha)$, $\forall u\ge-\alpha$. Hence
\begin{eqnarray*}
\E_{\Q^{(\alpha)}}\bigl(\xi_{n,E_n^c}^2\bigr) &\le&
{\Q^{(\alpha)}(E_n^c)
\over c_1^2 (1+\varepsilon n^{1/2}+\alpha)^2} + {1\over c_1 R_\alpha(0)} \E \biggl(
{\mathbf{1}_{\{ S_n < \varepsilon n^{1/2},
\underline{S}_n\ge-\alpha\} }
\over S_n+\alpha+1} \biggr)
\\
&=& o \biggl({1\over n} \biggr) + {1\over c_1 R_\alpha(0)} \E
\biggl( {\mathbf{1}_{\{ S_n < \varepsilon n^{1/2},
\underline{S}_n\ge-\alpha\} }
\over S_n+\alpha+1} \biggr),
\end{eqnarray*}
the last line following from the assumption that $\Q
^{(\alpha
)}(E_n^c) \to0$. For the expectation term on the right-hand side, we
observe that, by Lemma \ref{la1},
\begin{eqnarray*}
\E \biggl( {\mathbf{1}_{\{ S_n < \varepsilon n^{1/2},
\underline{S}_n\ge-\alpha\} }
\over S_n+\alpha+1} \biggr) &\le& \sum
_{i=0}^{\lceil\varepsilon n^{1/2} + \alpha\rceil-1} \E \biggl( {\mathbf{1}_{\{ -\alpha+ i \le S_n
< -\alpha+ i+1,
\underline{S}_n \ge-\alpha\} }
\over S_n+\alpha+1}
\biggr)
\\
&\le&\sum_{i=0}^{\lceil\varepsilon n^{1/2} + \alpha\rceil-1} {1\over i +1}
c_3 {(\alpha+1) (i+1)\over n^{3/2}}
\\
&=& {\lceil\varepsilon n^{1/2} + \alpha\rceil c_3 (\alpha+1)
\over n^{3/2}}.
\end{eqnarray*}

We have therefore proved that
\[
\E_{\Q^{(\alpha)}}\bigl(\xi_{n,E_n^c}^2\bigr) \le o \biggl(
{1\over n} \biggr) + {\lceil\varepsilon n^{1/2} + \alpha\rceil c_3 (\alpha+1) \over
n^{3/2} c_1 R_\alpha(0)},\qquad n\to\infty.
\]

 Since $\varepsilon$ can be arbitrarily small (whereas the
constants $c_1$ and $c_3$ do not depend on $\varepsilon$), this yields
Lemma \ref{lxi}.
\end{pf*}

The proof of Lemma \ref{lindependance} needs some preparation. We
start by the following elementary fact. Recall that $\log_+ y:= \max
\{
0, \log y\}$ for any $y\ge0$.

\begin{lemma}[(\cite{elie}, Lemma B.1)]
\label{lX-Xtilde}
Let $X \ge0$ and $\widetilde{X} \ge0$ be random variables such that
$\E[X \log_+^2 X] + \E[\widetilde{X} \log_+ \widetilde{X}]
<\infty$. Then
%
\begin{eqnarray}
\E \bigl[ X \log_+^2 \widetilde{X} \bigr] + \E [ \widetilde{X}
\log_+ X ] &<& \infty, \label{csq0-cond-X}
\\
\lim_{z\to\infty} {1\over z} \E \bigl[ X
\log_+^2 (X+\widetilde{X}) \min\bigl\{ \log_+ (X+\widetilde{X}), z
\bigr\} \bigr] &=& 0, \label{csq1-cond-X}
\\
\lim_{z\to\infty} {1\over z} \E \bigl[ \widetilde{X}
\log_+ (X+\widetilde{X}) \min\bigl\{ \log_+ (X+\widetilde{X}), z\bigr\} \bigr] &=&
0. \label{csq2-cond-X}
\end{eqnarray}
\end{lemma}

We continue our preparation for the proof of Lemma \ref
{lindependance}. Let $k_n< n$ be an integer such that $k_n \to\infty$
($n\to\infty$). Recall that we defined $W_n^{(\alpha)} = \sum_{|x|=n}
\ee^{-V(x)} \times \mathbf{1}_{\{\underline{V}(x)\ge-\alpha\} }$. For each
vertex $x$ with $|x|=n$ and $x\not= w^{(\alpha)}_n$, there is a unique
$i$ with $0\le i< n$ such that $w^{(\alpha)}_i \le x$ and that
$w^{(\alpha)}_{i+1} \not\le x$. For any $i\ge1$, let
\[
\Omega\bigl(w^{(\alpha)}_i\bigr):= \bigl\{ |x|= i \dvtx
x>w^{(\alpha)}_{i-1}, x\not= w^{(\alpha)}_i \bigr\}.
\]
[In words, $\Omega(w^{(\alpha)}_i)$ stands for the set of
``brothers'' of $w^{(\alpha)}_i$.] Accordingly,
\[
W_n^{(\alpha)} = \ee^{-V(w^{(\alpha)}_n)} \mathbf{1}_{\{\underline{V}(w^{(\alpha
)}_n)\ge
-\alpha\} } +
\sum_{i=0}^{n-1} \sum
_{y\in\Omega(w^{(\alpha)}_{i+1})} \sum_{|x|=n,
x\ge y}
\ee^{-V(x)} \mathbf{1}_{\{\underline{V}(x)\ge-\alpha\} }.
\]

 We write
\begin{eqnarray*}
W_n^{(\alpha), [0,k_n)} &:=& \sum_{i=0}^{k_n-1}
\sum_{y\in\Omega(w^{(\alpha)}_{i+1})} \sum_{|x|=n, x\ge y}
\ee^{-V(x)} \mathbf{1}_{\{\underline{V}(x)\ge-\alpha\} },
\\
W_n^{(\alpha), [k_n,n]} &:=& \ee^{-V(w^{(\alpha)}_n)} \mathbf{1}_{\{\underline{V}(w^{(\alpha)}_n)\ge-\alpha\} } +
\sum_{i=k_n}^{n-1} \sum
_{y\in\Omega(w^{(\alpha)}_{i+1})} \sum_{|x|=n, x\ge y}
\ee^{-V(x)} \mathbf{1}_{\{\underline{V}(x)\ge-\alpha\} },
\end{eqnarray*}
so that $W_n^{(\alpha)} = W_n^{(\alpha), [0,k_n)} +
W_n^{(\alpha), [k_n,n]}$. We define $D_n^{(\alpha), [0,k_n)}$ and
$D_n^{(\alpha), [k_n,n]}$ similarly. Let
\begin{eqnarray*}
E_{n,1} &:=& \bigl\{ k_n^{1/3} \le V
\bigl(w^{(\alpha)}_{k_n}\bigr) \le k_n\bigr\} \cap
\bigcap_{i=k_n}^n \bigl\{ V
\bigl(w^{(\alpha)}_i\bigr) \ge k_n^{1/6}\bigr
\},
\\
E_{n,2} &:=& \bigcap_{i=k_n}^{n-1}
\biggl\{ \sum_{y\in\Omega(w^{(\alpha)}_{i+1})} \bigl[1+\bigl(V(y)-V
\bigl(w^{(\alpha)}_i\bigr)\bigr)^+\bigr] \ee^{-[V(y)-V(w^{(\alpha)}_i)]} \le
\ee^{V(w^{(\alpha)}_i)/2} \biggr\},
\\
E_{n,3} &:=& \biggl\{ D_n^{(\alpha),[k_n,n]}\le
{1\over n^2} \biggr\}.
\end{eqnarray*}

We choose
%
\begin{equation}
E_n:= E_{n,1} \cap E_{n,2} \cap
E_{n,3}. \label{En}
\end{equation}

\begin{lemma}
\label{lQE->0}
Assume (\ref{cond-hab}),
(\ref{cond-2}) and (\ref{cond-X}). Let
$\alpha\ge0$. Let $k_n$ be such that ${k_n \over(\log n)^6}\to
\infty
$ and
that ${k_n\over n^{1/2}} \to0$, $n\to\infty$.
Let $E_n$ be as in (\ref{En}). Then
\[
\lim_{n\to\infty} \Q^{(\alpha)}(E_n)=1, \qquad\lim
_{n\to\infty} \inf_{u\in[k_n^{1/3}, k_n]} \Q^{(\alpha)}
\bigl(E_n | V\bigl(w^{(\alpha)}_{k_n}\bigr) =u\bigr) =1.
\]
\end{lemma}

\begin{pf} Write, for $i\ge0$,
\[
E_2^{(i)}:= \biggl\{ \sum_{y\in\Omega(w^{(\alpha)}_{i+1})}
\bigl[1+\bigl(V(y)-V\bigl(w^{(\alpha)}_i\bigr)\bigr)^+\bigr]
\ee^{-[V(y)-V(w^{(\alpha)}_i)]} \le\ee^{V(w^{(\alpha)}_i)/2} \biggr\}.
\]
(Thus $E_{n,2} = \bigcap_{i=k_n}^{n-1} E_2^{(i)}$.)

For $z\ge-\alpha$, let $\Q^{(\alpha)}_z$ be the law of $\mathcal
B_{\alpha}$ (in Fact \ref{pspine}) when the ancestor particle is
located at position $z$. (So $\Q^{(\alpha)}_0 = \Q^{(\alpha)}$.) We
claim that
%
\begin{eqnarray}
\sum_{i\ge0} \Q^{(\alpha)}_z \bigl[
\bigl(E_2^{(i)}\bigr)^c\bigr] &<& \infty\qquad \forall
z\ge-\alpha, \label{eqborel}
\\
\lim_{z\to\infty} \sum_{i\ge0}
\Q^{(\alpha)}_z \bigl[\bigl(E_2^{(i)}
\bigr)^c\bigr] &=& 0. \label{eqborelz}
\end{eqnarray}

To check (\ref{eqborel}) and (\ref{eqborelz}), we observe that by
Fact \ref{pspine}, for any integer $i\ge0$ and real number $u\ge
-\alpha$,
\begin{eqnarray*}
&&\Q^{(\alpha)}_z \bigl[\bigl(E_2^{(i)}
\bigr)^c | V\bigl(w^{(\alpha)}_i\bigr) = u\bigr] \\[-2pt]
&&\qquad=
\Q^{(\alpha)}_u \biggl\{ \sum_{x\in\Omega(w^{(\alpha)}_1)}
\bigl[1+ \bigl(V(x)-u\bigr)^+\bigr]\ee^{-[V(x)-u]} >\ee^{u/2} \biggr\}
\\[-2pt]
&&\qquad\le \Q^{(\alpha)}_u \biggl\{ \sum
_{|x|=1} \bigl[1+ \bigl(V(x)-u\bigr)^+\bigr]\ee ^{-[V(x)-u]} >
\ee^{u/2} \biggr\}.
\end{eqnarray*}

So, if $\E_u$ denotes expectation with respect to the law of
the branching random walk with the ancestor particle located at $u$, then
\begin{eqnarray*}
&&\Q^{(\alpha)}_z \bigl[\bigl(E_2^{(i)}
\bigr)^c | V\bigl(w^{(\alpha)}_i\bigr) = u\bigr]
\\[-2pt]
&&\qquad\le \E_u \biggl[ {\sum_{|y|=1} R_\alpha(V(y)) \ee^{-V(y)} \mathbf{
1}_{\{
V(y) \ge-\alpha\} } \over R_\alpha(u) \ee^{-u}} \\[-2pt]
&&\hspace*{48pt}{}\times\mathbf{1}_{\{ \sum
_{|x|=1} [1+ (V(x)-u)^+]\ee^{-[V(x)-u]} >\ee^{u/2} \} }
\biggr]
\\[-2pt]
&&\qquad=\E \biggl[ {\sum_{|y|=1} R_\alpha(V(y)+u) \ee^{-[V(y)+u]} \mathbf{
1}_{\{
V(y) \ge-\alpha-u\} } \over R_\alpha(u) \ee^{-u}} \\[-2pt]
&&\hspace*{116pt}{}\times \mathbf{1}_{\{ \sum
_{|x|=1} [1+ V(x)^+]\ee^{- V(x)} >\ee^{u/2} \} } \biggr].
\end{eqnarray*}

By (\ref{h0>}), there exists a constant $c_{10}>0$ such that
\[
{R_\alpha(V(y)+u) \over R_\alpha(u)} \le c_{10} {V(y)^+ + u + \alpha+
1\over u+\alpha+1} = c_{10} \biggl[1+ {V(y)^+\over u+\alpha+1}\biggr];
\]
thus
\begin{eqnarray*}
&&\Q^{(\alpha)}_z \bigl[\bigl(E_2^{(i)}
\bigr)^c | V\bigl(w^{(\alpha)}_i\bigr) = u\bigr]\\
&&\qquad\le
c_{10} \E \biggl[ \sum_{|y|=1}
\ee^{-V(y)} \mathbf{1}_{\{ \sum
_{|x|=1} [1+ V(x)^+]\ee^{- V(x)} >\ee^{u/2} \} }
\\
&&\hspace*{27pt}\qquad\quad{} + {1\over u+\alpha+1} \sum_{|y|=1} V(y)^+
\ee^{-V(y)} \mathbf{1}_{\{
\sum_{|x|=1} [1+ V(x)^+]\ee^{- V(x)} >\ee^{u/2} \} } \biggr]
\\
&&\qquad= c_{10} \E \biggl[ X \mathbf{1}_{\{ X + \widetilde{X} >\ee^{u/2} \} } +
{\widetilde{X} \mathbf{1}_{\{ X + \widetilde{X} >\ee^{u/2} \} }\over
u+\alpha+1} \biggr],
\end{eqnarray*}
where $X:= \sum_{|y|=1} \ee^{-V(y)}$ and $\widetilde{X}:=
\sum_{|y|=1} V(y)^+ \ee^{-V(y)}$. Consequently,
\[
\Q^{(\alpha)}_z \bigl[\bigl(E_2^{(i)}
\bigr)^c\bigr] \le c_{10} \bigl(\E\otimes
\E^{(\alpha)}_z\bigr) \biggl[ X \mathbf{1}_{\{ X +
\widetilde
{X} >\ee^{S_i/2} \} } +
{\widetilde{X} \mathbf{1}_{\{ X + \widetilde{X}
>\ee^{S_i/2} \} }\over S_i+\alpha+1} \biggr],
\]
where, on the right-hand side, we assume that $(X, \widetilde
{X})$ and $S_i$ are independent, the expectation $\E$ being for $(X,
\widetilde{X})$, while the expectation $\E^{(\alpha)}_z$ for $S_i$.
Here, $\E^{(\alpha)}_z$ stands for the expectation with respect to
$\P
^{(\alpha)}_z$, the law of the $h$-process of $(S_i)$ starting from $z$
and conditioned to stay in $[-\alpha, \infty)$; the transition
probabilities of this $h$-process being given in (\ref{h-process}).

Let us consider the expression on the right-hand side. We first take
the expectation for $S_i$ with respect to $\E^{(\alpha)}_z$.
The event $\{ X + \widetilde{X} >\ee^{S_i/2} \}$ can be written as $S_i
< 2 \log(X + \widetilde{X})$. Therefore, by the definition of $\E
^{(\alpha)}_z$, for any $x\ge0$ and $\widetilde{x}\ge0$,
\begin{eqnarray*}
&&\E^{(\alpha)}_z \biggl[ x \mathbf{1}_{\{ x + \widetilde{x} >\ee^{S_i/2}
\} } +
{\widetilde{x} \mathbf{1}_{\{ x + \widetilde{x} >\ee^{S_i/2} \}
}\over S_i+\alpha+1} \biggr]
\\
&&\qquad= {1\over R_\alpha(z)} \E \biggl[ R_\alpha(S_i+z) {
\mathbf{1}}_{\{ \underline{S}_i \ge-z-\alpha
\} }\\
&&\hspace*{72pt}{}\times \biggl( x \mathbf{1}_{\{ S_i+z < 2 \log(x + \widetilde{x}) \} } + {\widetilde{x} \mathbf{1}_{\{ S_i+z < 2 \log(x + \widetilde{x}) \}
}\over
S_i + z +\alpha+1}
\biggr) \biggr],
\end{eqnarray*}
which, by (\ref{c0}), is\footnote{The constant $c_{11}$, as
well as the forthcoming $c_{12}$ and $c_{13}$, may depend on $\alpha$.
This, however, makes no trouble as $\alpha$ will ultimately be a large
(but fixed) constant.}
\begin{eqnarray*}
&\le& {c_2 \over R_\alpha(z)} \E \biggl[ (S_i + z +\alpha+1) {
\mathbf{1}}_{\{ \underline{S}_i \ge
-z-\alpha
\} }\\
&&\hspace*{19pt}\qquad{}\times \biggl( x \mathbf{1}_{\{ S_i+z< 2 \log(x + \widetilde{x}) \} } + {\widetilde{x} \mathbf{1}_{\{ S_i +z< 2 \log(x + \widetilde{x}) \}
}\over
S_i + z +\alpha+1}
\biggr) \biggr]
\\
&\le&{c_{11} [x(1+\log_+ (x + \widetilde{x})) + \widetilde{x}]
\over
R_\alpha(z)} \P \bigl\{ \underline{S}_i \ge-z-
\alpha, S_i +z < 2 \log(x + \widetilde {x}) \bigr\}.
\end{eqnarray*}

 Applying Lemma \ref{la6} yields that
\begin{eqnarray*}
\hspace*{-4pt}&&\sum_{i\ge0} \E^{(\alpha)}_z
\biggl[ x \mathbf{1}_{\{ x + \widetilde{x}
>\ee^{S_i/2} \} } + {\widetilde{x} \mathbf{1}_{\{ x + \widetilde{x}
>\ee
^{S_i/2} \} }\over S_i+\alpha+1} \biggr]
\\
\hspace*{-9pt}&&\qquad\le\hspace*{-1pt}{c_{12} [x(1+\log_+ (x + \widetilde{x})) + \widetilde{x}] [1+
\log
_+ (x + \widetilde{x})] [1+ \min\{ \log_+ (x + \widetilde{x}), z\}]
\over R_\alpha(z)}.
\end{eqnarray*}

Taking expectation for $(X, \widetilde{X})$, using (\ref
{csq0-cond-X})--(\ref{csq2-cond-X}) in Lemma \ref{lX-Xtilde} [which we
are entitled to apply, in view of assumption (\ref{cond-X})], and
recalling from (\ref{c0}) that $R_\alpha(z)$ grows linearly when
$z\to
\infty$, we obtain (\ref{eqborel}) and (\ref{eqborelz}).

We now prove that $\Q^{(\alpha)}(E_n) \to1$, $n\to\infty$. Since $E_n
= E_{n,1} \cap E_{n,2} \cap E_{n,3}$, let us check that $\lim_{n\to
\infty} \Q^{(\alpha)}(E_{n,\ell})=1$, for $\ell= 1$ and $2$, and that
$\lim_{n\to\infty} \Q^{(\alpha)}(E_{n,3}^c \cap E_{n,1}\cap E_{n,2})=0$.

For $E_{n,1}$: Fact \ref{pspine} says that $(V(w^{(\alpha)}_n), n\ge
0)$ under $\Q^{(\alpha)}$ is the centered random walk $(S_n)$
conditioned to stay in $[-\alpha, \infty)$; so it is clear that $\Q
^{(\alpha)}(E_{n,1}) \to1$, $n\to\infty$.

For $E_{n,2}$: this follows from (\ref{eqborel}) (by taking $z=0$ there).

For $E_{n,3}$: Let $\mathcal G_\infty:=\sigma\{ V(w_k^{(\alpha)}),
V(z), z\in\Omega(w^{(\alpha)}_{k+1}), k\ge0\}$ be the sigma-algebra
generated by the positions of the spine and its brothers. We know that
the branching random walk rooted at $z\in\Omega(w^{(\alpha)}_i)$ has
the same law under $\P$ and under $\Q^{(\alpha)}$. Therefore,
\[
\E_{\Q^{(\alpha)}}\bigl[D_n^{(\alpha),[k_n,n]} | \mathcal
G_\infty\bigr] = R_\alpha\bigl(V\bigl(w_n^{(\alpha)}
\bigr)\bigr)\ee^{-V(w_n^{(\alpha)})} + \sum_{i=k_n}^{n-1}
\sum_{z\in\Omega(w^{(\alpha)}_{i+1})} R_\alpha \bigl(V(z)\bigr)
\ee^{-V(z)}.
\]

 For $z\in\Omega(w^{(\alpha)}_{i+1})$, we have $R_\alpha
(V(z)) \le c_{13} [1+\alpha+V(w_i^{(\alpha)})] [
1+(V(z)-V(w_i^{(\alpha
)}))^+]$. Therefore,
%
\begin{equation}
\mathbf{1}_{E_{n,1}\cap E_{n,2}} \E_{\Q^{(\alpha)}}\bigl[D_n^{(\alpha),[k_n,n]} |
\mathcal G_\infty\bigr] = O \bigl(n \ee^{-k_n^{1/6}/3} \bigr),\qquad n\to\infty,
\label{eqEn3}
\end{equation}
where the $O(n \ee^{-k_n^{1/6}/3})$ term on the right-hand
side represents a deterministic expression. Since ${k_n\over(\log
n)^6} \to\infty$, it follows from the Markov inequality that $\Q
^{(\alpha)}(E_{n,3}^c \cap E_{n,1}\cap E_{n,2}) \to0$, $n\to\infty$.

It remains to check that $\Q^{(\alpha)}(E_n | V(w^{(\alpha)}_{k_n}) =u)
\to1$ uniformly in $u\in[k_n^{1/3}, k_n]$.

By (\ref{eqborelz}), $\Q^{\alpha}(E_{n,2}^c | V(w^{(\alpha)}_{k_n})
=u) \to0$ uniformly in $u\in[k_n^{1/3}, k_n]$, whereas according to
(\ref{eqEn3}), $\mathbf{1}_{E_{n,1}\cap E_{n,2}} \Q^{(\alpha)}(E_{n,3}^c
| \mathcal G_\infty)$ is bounded by a deterministic expression which
goes to 0 when $n\to\infty$. Therefore, we only have to check that
$\Q
^{(\alpha)} ( E_{n,1} | V(w^{(\alpha)}_{k_n}) = u ) \to1$, uniformly
in $u\in[k_n^{1/3}, k_n]$. By Fact \ref{pspine} and (\ref{h-process}),
\[
\Q^{(\alpha)} \bigl( E_{n,1} | V\bigl(w^{(\alpha)}_{k_n}
\bigr) = u \bigr) = {1\over R_\alpha(u)} \E\bigl[R_\alpha(S_{n-k_n}+u)
\mathbf{1}_{\{ \underline{S}_{n-k_n} \ge k_n^{1/6} -u\}}\bigr].
\]

 Let, as before, $c_0:=\lim_{t\to\infty} {R_\alpha(t)\over
t}$, and let $\eta\in(0, c_0)$. Let $f_\eta(t):= (c_0 - \eta) \min
\{ t, {1\over\eta}\}$. Then $R_\alpha(t) \ge b f_\eta({t\over b})$
for all sufficiently large $t$ and uniformly in $b>0$. We take $b:=
(n-k_n)^{1/2}\sigma$ (with $\sigma^2:= \E[S_1^2]$ as before), to see
that for all sufficiently large $n$ and uniformly in $u>k_n^{1/6}$,
\begin{eqnarray*}
&&\Q^{(\alpha)} \bigl( E_{n,1} | V\bigl(w^{(\alpha)}_{k_n}
\bigr) = u \bigr) \\
&&\qquad\ge {(n-k_n)^{1/2} \sigma\over R_\alpha(u)} \E \biggl[ f_\eta
\biggl( {S_{n-k_n}+u \over(n-k_n)^{1/2} \sigma} \biggr) \mathbf{1}_{\{ \underline{S}_{n-k_n} \ge k_n^{1/6} -u\} } \biggr]
\\
&&\qquad\ge {(n-k_n)^{1/2} \sigma\over R_\alpha(u)} \E \biggl[ f_\eta \biggl(
{S_{n-k_n}+u-k_n^{1/6} \over(n-k_n)^{1/2} \sigma} \biggr) \mathbf{1}_{\{ \underline{S}_{n-k_n} \ge k_n^{1/6} -u\}} \biggr].
\end{eqnarray*}

 Since ${k_n\over n^{1/2}}\to0$, we can apply (\ref
{elie-bruno}) to see that, as $n\to\infty$,
\[
\E \biggl[ f_\eta \biggl( {S_{n-k_n}+u-k_n^{1/6} \over
(n-k_n)^{1/2}\sigma
} \biggr) {
\mathbf{1}}_{\{ \underline{S}_{n-k_n} \ge k_n^{1/6} -u\}} \biggr] \sim {\theta R(u-k_n^{1/6}) \over(n-k_n)^{1/2}}\int
_0^\infty t \ee^{-t^2/2} f_\eta(t)
\,\d t,
\]
uniformly in $u\in[k_n^{1/6}, k_n]$. Consequently,
\[
\liminf_{n\to\infty} \inf_{u\in[k_n^{1/3}, k_n]} \Q^{(\alpha)}
\bigl( E_{n,1} | V\bigl(w^{(\alpha)}_{k_n}\bigr) = u
\bigr) \ge \theta\sigma\int_0^\infty t
\ee^{-t^2/2} f_\eta(t)\, \d t.
\]
Note that $\int_0^\infty t \ee^{-t^2/2} f_\eta(t) \,\d t \ge
(c_0-\eta) \int_0^{1/\eta}t^2 \ee^{-t^2/2} \,\d t$. Letting $\eta\to
0$ gives
\[
\liminf_{n\to\infty} \inf_{u\in[k_n^{1/3}, k_n]} \Q^{(\alpha)}
\bigl( E_{n,1} | V\bigl(w^{(\alpha)}_{k_n}\bigr) = u
\bigr) \ge c_0\theta\sigma \biggl({\pi\over2}
\biggr)^{ 1/2} =1,
\]
the last identity following from (\ref{relation-constants}).
Consequently, $\Q^{(\alpha)} ( E_n | V(w^{(\alpha)}_{k_n}) = u ) \to1$
uniformly in $u\in[k_n^{1/3}, k_n]$. Lemma \ref{lQE->0} is
proved.
\end{pf}

We now proceed to prove Lemma \ref{lindependance}.

\begin{pf*}{Proof of Lemma \ref{lindependance}} Let $k_n$ be such
that $k_n \to\infty$ and that ${k_n\over n^{1/2}} \to0$, $n\to
\infty
$. Let $E_n$ be the event in (\ref{En}). By Lemma \ref{lQE->0},
$\Q
^{(\alpha)}(E_n) \to1$, $n\to\infty$.

On $E_n$, we have $D_n^{(\alpha), [k_n,n]} \le{1\over n^2}$; in
particular, since $W_n^{(\alpha), [k_n,n]}\le D_n^{(\alpha), [k_n,n]}$,
we have $W_n^{(\alpha), [k_n,n]}\le{1\over n^2}$ on $E_n$. On the
other hand, $R_\alpha(V(w^{(\alpha)}_n)) \ge1$, so
%
\begin{eqnarray}\label{EWD1}
&&\E_{\Q^{(\alpha)}} \biggl[ {W_n^{(\alpha), [k_n,n]}\over D_n^{(\alpha)}} {\mathbf{1}_{E_n} \over R_\alpha(V(w^{(\alpha)}_n))}
\biggr]
\nonumber
\\[-8pt]
\\[-8pt]
\nonumber
&&\qquad \le \E_{\Q^{(\alpha)}} \biggl[ {{1/ n^2}\over D_n^{(\alpha)}} \biggr] = \E
\biggl[ {{1/ n^2}\over R_\alpha(0)} \biggr] = o \biggl( {1\over n}
\biggr).
\end{eqnarray}

It remains to treat ${W_n^{(\alpha), [0,k_n)}\over D_n^{(\alpha)}}
{\mathbf{1}_{E_n} \over R_\alpha(V(w^{(\alpha)}_n))}$. Since
$D_n^{(\alpha
)} \ge D_n^{(\alpha), [0,k_n)}$, we have\footnote{Notation: ${0\over0}:=0$ for the ratio ${W_n^{(\alpha), [0,k_n)}\over D_n^{(\alpha),
[0,k_n)}}$; noting that if $D_n^{(\alpha), [0,k_n)}=0$, then
$W_n^{(\alpha), [0,k_n)}=0$.}
\[
\E_{\Q^{(\alpha)}} \biggl[ {W_n^{(\alpha), [0,k_n)}
\over D_n^{(\alpha)}} {\mathbf{1}_{E_n} \over R_\alpha(V(w^{(\alpha)}_n))}
\biggr] \le \E_{\Q^{(\alpha)}} \biggl[ {W_n^{(\alpha), [0,k_n)}\over
D_n^{(\alpha), [0,k_n)}}
{\mathbf{1}_{E_n} \over R_\alpha(V(w^{(\alpha)}_n))} \biggr].
\]

 Therefore, by Fact \ref{pspine},
%
\begin{eqnarray}\label{EWD2}\quad
&&\E_{\Q^{(\alpha)}} \biggl[ {W_n^{(\alpha), [0,k_n)}
\over D_n^{(\alpha)}} {\mathbf{1}_{E_n} \over R_\alpha(V(w^{(\alpha)}_n))}
\biggr]
\nonumber
\\[-4pt]
\\[-12pt]
\nonumber
&&\qquad\le \E_{\Q^{(\alpha)}} \biggl( {W_n^{(\alpha), [0,k_n)}
\over D_n^{(\alpha), [0,k_n)}} {
\mathbf{1}}_{\{V(w^{(\alpha)}_{k_n}) \in[ k_n^{1/3}, k_n]\}} \biggr) \sup_{u\in[k_n^{1/3}, k_n]} \E^{(\alpha)}_u
\biggl( {1\over R_\alpha(S_{n-k_n})} \biggr).
\end{eqnarray}

For any $u\ge-\alpha$ and $j\ge1$, we have $\E^{(\alpha)}_u (
{1\over
R_\alpha(S_j)} ) = {1\over R_\alpha(u)} \P\{ \underline{S}_j \ge
-\alpha
-u\}$, which yields, by (\ref{kozlov-uniform}),
\[
\sup_{u\in[k_n^{1/3}, k_n]} \E^{(\alpha)}_u \biggl(
{1\over R_\alpha
(S_{n-k_n})} \biggr) \sim{\theta\over(n-k_n)^{1/2}} \sim
{\theta
\over
n^{1/2}},\qquad n\to\infty.
\]

Going back to (\ref{EWD2}), we obtain
\begin{eqnarray*}
&&\E_{\Q^{(\alpha)}} \biggl[ {W_n^{(\alpha), [0,k_n)}
\over D_n^{(\alpha)}} {\mathbf{1}_{E_n} \over R_\alpha(V(w^{(\alpha)}_n))}
\biggr] \\
&&\qquad\le {\theta+ o(1) \over n^{1/2}} \E_{\Q^{(\alpha)}} \biggl(
{W_n^{(\alpha), [0,k_n)}
\over D_n^{(\alpha), [0,k_n)}} \mathbf{1}_{\{V(w^{(\alpha)}_{k_n})\in[k_n^{1/3}, k_n]\}} \biggr).
\end{eqnarray*}

 We claim that
%
\begin{equation}
\limsup_{n\to\infty} n^{1/2} \E_{\Q^{(\alpha)}} \biggl(
{W_n^{(\alpha), [0,k_n)}
\over D_n^{(\alpha), [0,k_n)}} \mathbf{1}_{\{V(w^{(\alpha)}_{k_n}) \in[k_n^{1/3}, k_n]\}} \biggr) \le \theta. \label{claim}
\end{equation}

 Then we will have
\[
\E_{\Q^{(\alpha)}} \biggl[ {W_n^{(\alpha), [0,k_n)}
\over D_n^{(\alpha)}} {\mathbf{1}_{E_n} \over R_\alpha(V(w^{(\alpha)}_n))}
\biggr] \le {\theta^2\over n} + o \biggl( {1\over n} \biggr),
\]
which, together with (\ref{EWD1}) and remembering
$W_n^{(\alpha)} = W_n^{(\alpha), [0,k_n)} + W_n^{(\alpha), [k_n,n]}$,
will complete the proof of Lemma \ref{lindependance}.

It remains to check (\ref{claim}). By Fact \ref{pspine},
\begin{eqnarray*}
&&\E_{\Q^{(\alpha)}} \biggl( {W_n^{(\alpha), [0,k_n)}
\over D_n^{(\alpha), [0,k_n)}} \mathbf{1}_{E_n}
\biggr)
\\
&&\quad\ge \E_{\Q^{(\alpha)}} \biggl( {W_n^{(\alpha), [0,k_n)}
\over D_n^{(\alpha), [0,k_n)}} {
\mathbf{1}}_{\{V(w^{(\alpha)}_{k_n}) \in[k_n^{1/3}, k_n]\}} \biggr) \inf_{u\in[k_n^{1/3}, k_n]} \Q^{(\alpha)}
\bigl( E_n | V\bigl(w^{(\alpha)}_{k_n}\bigr) = u \bigr).
\end{eqnarray*}

 By Lemma \ref{lQE->0}, $\inf_{u\in[k_n^{1/3}, k_n]}\Q
^{(\alpha)} ( E_n | V(w^{(\alpha)}_{k_n}) = u ) \to1$. Therefore,
as\break
$n\to\infty$,
\[
\E_{\Q^{(\alpha)}} \biggl( {W_n^{(\alpha), [0,k_n)}
\over D_n^{(\alpha), [0,k_n)}} \mathbf{1}_{\{V(w^{(\alpha)}_{k_n}) \in[k_n^{1/3}, k_n]\}}
\biggr) \le \bigl(1+ o(1)\bigr) \E_{\Q^{(\alpha)}} \biggl( {W_n^{(\alpha), [0,k_n)}
\over D_n^{(\alpha), [0,k_n)}}
\mathbf{1}_{E_n} \biggr).
\]

 Since $D_n^{(\alpha), [0,k_n)}\ge W_n^{(\alpha), [0,k_n)}$,
we have
\begin{eqnarray*}
&&\E_{\Q^{(\alpha)}} \biggl( {W_n^{(\alpha), [0,k_n)}
\over D_n^{(\alpha), [0,k_n)}} \mathbf{1}_{E_n}
\biggr)\\
&&\qquad \le \E_{\Q^{(\alpha)}} \biggl( {W_n^{(\alpha), [0,k_n)}
\over D_n^{(\alpha), [0,k_n)}} {
\mathbf{1}}_{E_n}\mathbf{1}_{\{D_n^{(\alpha)}>
{1\over n}\}} \biggr) + \Q^{(\alpha)} \biggl(
D_n^{(\alpha)}\le{1\over n} \biggr).
\end{eqnarray*}

 Let $0<\eta_1<1$. By the Markov inequality, we see that\break  $\Q
^{(\alpha)} (D_n^{(\alpha)}\le{1\over n}) \le{1\over n} \E_{\Q
^{(\alpha)}}({1\over D_n^{(\alpha)}}) = {1\over n R_\alpha(0)}$. On the
other hand, we already noticed that $D_n^{(\alpha), [k_n,n]} \mathbf{
1}_{E_n}$ is bounded by a deterministic $o({1 \over n})$. Therefore,
for all sufficiently large $n$, $D_n^{(\alpha), [k_n,n]} \le\eta_1
D_n^{(\alpha)}$ on $E_n\cap\{D_n^{(\alpha)}>{1\over n}\}$.
Accordingly, for all sufficiently large~$n$,
\begin{eqnarray*}
\E_{\Q^{(\alpha)}} \biggl( {W_n^{(\alpha), [0,k_n)}
\over D_n^{(\alpha), [0,k_n)}} \mathbf{1}_{E_n}
\biggr) &\le& {1\over1-\eta_1} \E_{\Q^{(\alpha)}} \biggl(
{W_n^{(\alpha), [0,k_n)}
\over D_n^{(\alpha)}} \mathbf{1}_{E_n \cap\{D_n^{(\alpha)}> {1\over n}\}
} \biggr) + {1\over n R_\alpha(0)}
\\
&\le& {1\over1-\eta_1} \E_{\Q^{(\alpha)}} \biggl( {W_n^{(\alpha)}
\over D_n^{(\alpha)}}
\biggr) + {1\over n R_\alpha(0)}.
\end{eqnarray*}

 On the right-hand side, $\E_{\Q^{(\alpha)}}
( {W_n^{(\alpha)} \over D_n^{(\alpha)}} ) \sim{\theta\over n^{1/2}}$;
see (\ref{moment1}). It follows that
\[
\limsup_{n\to\infty} n^{1/2} \E_{\Q^{(\alpha)}} \biggl(
{W_n^{(\alpha), [0,k_n)}
\over D_n^{(\alpha), [0,k_n)}} \mathbf{1}_{\{V(w^{(\alpha)}_{k_n}) \in[k_n^{1/3}, k_n]\}} \biggr) \le{\theta\over1-\eta_1}.
\]

Sending $\eta_1 \to0$ gives (\ref{claim}), and completes the
proof of Lemma \ref{lindependance}.
\end{pf*}

\begin{pf*}{Proof of Proposition \ref{pmoments}} Equation
(\ref{moment2}) follows from Lemmas \ref{lxi} and~\ref
{lindependance}.
\end{pf*}

\section{\texorpdfstring{Proof of Theorem \protect\ref{tmain}}
{Proof of Theorem 1.1}}
\label{sproof}

Assume (\ref{cond-hab}), (\ref{cond-2}) and (\ref
{cond-X}). Let $\alpha\ge0$. By Proposition \ref{pmoments}, under
$\Q
^{(\alpha)}$, $n^{1/2}{W_n^{(\alpha)}\over D_n^{(\alpha)}}$ converges,
as $n\to\infty$, in probability to $\theta$. Therefore, for any
$0<\varepsilon<1$,
\[
\Q^{(\alpha)} \biggl\{ \biggl| n^{1/2}{W_n^{(\alpha)}\over D_n^{(\alpha
)}} - \theta\biggr |
> \theta\varepsilon \biggr\} \to0,\qquad n\to\infty,
\]
that is,
\[
\E \bigl[ D_n^{(\alpha)} \mathbf{1}_{ \{ | n^{1/2}({W_n^{(\alpha)}/
D_n^{(\alpha)}}) - \theta| > \theta\varepsilon\} } \bigr] \to0,\qquad n
\to \infty.
\]

Recall that $\P^* (\bullet):= \P( \bullet| \mbox{nonextinction})$. By
Biggins~\cite{biggins-lindley}, condition\break  $\E(
\sum_{|x|=1} \ee^{-V(x)} ) =1$ in (\ref{cond-hab}) implies that
$\inf_{|x|=n} V(x) \to\infty$, $\P^*$-a.s.; thus $\inf_{|x|\ge0} V(x)
>-\infty$, $\P^*$-a.s.

Let $\Omega_k:= \{ \inf_{|x| \ge0} V(x) \ge-k\} \cap\{ \mbox
{nonextinction}\}$. Then $(\Omega_k, k\ge1)$ is a sequence of
nondecreasing events such that $\P^*(\bigcup_{k\ge1} \Omega_k)= \P
^*(\mbox{nonextinction}) = 1$. Let $\eta>0$. There exists
$k_0=k_0(\eta
)$ such that $\P^*( \Omega_{k_0}) \ge1-\eta$.

Since $\mathbf{1}_{\Omega_{k_0}}\le1$, we have
\[
\E \bigl[ D_n^{(\alpha)} \mathbf{1}_{ \{ | n^{1/2}({W_n^{(\alpha)}/
D_n^{(\alpha)}} )- \theta| > \theta\varepsilon\} } {
\mathbf{1}}_{\Omega
_{k_0}} \bigr] \to0,\qquad n\to\infty.
\]

 Because $D_n^{(\alpha)} \ge0$, this is equivalent to say
that, under $\P$,
%
\begin{eqnarray}\label{L1}
D_n^{(\alpha)} \mathbf{1}_{ \{
| n^{1/2}({W_n^{(\alpha)}/ D_n^{(\alpha)}}) - \theta| >
\theta\varepsilon\} } \mathbf{1}_{\Omega_{k_0}}
\to0
\nonumber
\\[-8pt]
\\[-8pt]
\eqntext{\mbox{in $L^1(\P)$, a fortiori in probability}.}
\end{eqnarray}

On $\Omega_{k_0}$, we have $W_n^{(\alpha)}=W_n$ for all $n$ and all
$\alpha\ge k_0$. For the behavior of $D_n^{(\alpha)}$, we observe that
according to (\ref{c0}), there exists a constant $M=M(\varepsilon)>0$
sufficiently large such that
\[
c_0 (1-\varepsilon) u \le R(u) \le c_0(1+\varepsilon)
u\qquad \forall u\ge M.
\]

 We fix our choice of $\alpha$ from now on: $\alpha:= k_0 +
M$. Since $R_\alpha(u) = R(u+\alpha)$, we have, on $\Omega_{k_0}$,
$0<c_0 (1-\varepsilon)(V(x)+\alpha) \le R_\alpha(V(x)) \le c_0
(1+\varepsilon)(V(x)+\alpha)$ (for all vertices $x$), so that on
$\Omega_{k_0}$,
\[
0<c_0 (1-\varepsilon) (D_n + \alpha W_n)
\le D_n^{(\alpha)} \le c_0 (1+\varepsilon)
(D_n + \alpha W_n)  \qquad\forall n.
\]

(We insist on the fact that on $\Omega_{k_0}$, $D_n + \alpha
W_n>0$ for all $n$.)

Recall that $D_n \to\mathscr{W}^*>0$, $\P^*$-a.s., and that $W_n \to
0$, $\P^*$-a.s. Therefore, on the one hand, $\liminf_{n\to\infty}
D_n^{(\alpha)} \ge c_0 (1-\varepsilon)\mathscr{W}^*>0$, $\P
^*$-a.s. on
$\Omega_{k_0}$; on the other hand, on $\Omega_{k_0}$,
\[
A_n \subset \biggl\{ \biggl| n^{1/2}{W_n^{(\alpha)} \over D_n^{(\alpha)}} -
\theta\biggr| > \theta\varepsilon \biggr\} \qquad \forall n,
\]
where
\[
A_n:= \biggl\{ n^{1/2} {W_n \over D_n + \alpha W_n} > (1+
\varepsilon)^2 c_0 \theta \biggr\} \cup \biggl\{
n^{1/2} {W_n \over D_n + \alpha W_n} < (1-\varepsilon)^2
c_0 \theta \biggr\}.
\]

 In view of (\ref{L1}), we obtain that, under $\P^*$,
\[
\mathbf{1}_{A_n} \mathbf{1}_{\Omega_{k_0}} \to0 \qquad\mbox{in probability},
\]
that is, $\P^* (A_n \cap\Omega_{k_0}) \to0$, $n\to\infty$.
Since $\P^* (\Omega_{k_0}) \ge1- \eta$, this implies
\[
\limsup_{n\to\infty} \P^* (A_n) \le\eta.
\]

 In other words, $n^{1/2}{W_n \over D_n}$ converges in
probability (under $\P^*$) to $c_0 \theta$, which is $({2\over\pi
\sigma^2})^{1/2}$ according to (\ref{relation-constants}). Theorem
\ref
{tmain} now follows by an application of Theorem \ref{thB} in the
\hyperref[s:intro]{Introduction}.

\section{\texorpdfstring{Proof of Theorem \protect\ref{tmain2}}
{Proof of Theorem 1.2}}
\label{snon-as}

We first study the minimal displacement in a branching
random walk. Recall that $\P^* (\bullet):= \P( \bullet| \mbox{nonextinction})$.

\begin{theorem}
\label{tminimal-position}
Assume (\ref{cond-hab}), (\ref{cond-2}) and (\ref{cond-X}).
We have
\[
\liminf_{n\to\infty} \biggl( \min_{|x|=n} V(x) -
{1\over2}\log n \biggr) = -\infty,\qquad \mbox{$\P^*$-a.s.}
\]
\end{theorem}

\begin{remark*}Although we are not going to use it, we mention
that\break  $\min_{|x|=n} V(x)$ behaves typically like ${3\over2}\log n$: if
conditions (\ref{cond-hab}), (\ref{cond-2}) and (\ref{cond-X}) hold,
then under $\P^*$, ${1\over\log n} \min_{|x|=n} V(x)\to{3\over2}$ in
probability; see \cite{yzpolymer}, \cite{addario-berry-reed} or \cite
{ezsimple} for proofs under some additional assumptions. A proof
assuming only (\ref{cond-hab}), (\ref{cond-2}) and (\ref{cond-X}) can
be found in \cite{elie}. In particular, we cannot replace ``$\liminf$''
in Theorem \ref{tminimal-position} by ``$\lim$.''
\end{remark*}

By admitting Theorem \ref{tminimal-position} for the time being, we
are ready to prove Theorem~\ref{tmain2}.

\begin{pf*}{Proof of Theorem \ref{tmain2}} By definition,
$W_n =
\sum_{|x|=n} \ee^{-V(x)} \ge\break  \exp[-\min_{|x|=n} V(x)]$, so Theorem
\ref
{tmain2} is a consequence of Theorem \ref{tminimal-position}.
\end{pf*}

The rest of the section is devoted to the proof of Theorem \ref
{tminimal-position}. We use once again a change-of-probabilities
technique. This time, however, we only need the well-known
change-of-probabilities setting in Lyons~\cite{lyons}: Under (\ref
{cond-hab}), $(W_n)$ is a nonnegative martingale, so we can define a
probability $\Q$ such that for any $n$,
%
\begin{equation}
\Q|_{\mathscr{F}_n}:= W_n \bullet\P|_{\mathscr{F}_n}. \label{Q}
\end{equation}

 Recall that the positions of the particles in the first
generation, $(V(x),\break  |x|=1)$, are distributed under $\P$ as the point
process $\Theta$; let $\widehat{\Theta}$ denote a point process whose
distribution is the law of $(V(x), |x|=1)$ under $\Q$.

Lyons's spinal decomposition describes the distribution of the
branching random walk under $\Q$; it involves a spine process denoted
by $(w_n, n\ge0)$: We take $w_0:= \varnothing$, and the system starts
at the initial position $V(w_0)=0$. At time $1$, $w_0$~gives birth to
the point process $\widehat{\Theta}$. We choose $w_1$ at step $1$ among
the offspring~$x$ with probability proportional to $\ee^{-V(x)}$. The
particle $w_1$ gives birth to particles distributed as $\widehat
{\Theta
}$ [with respect to their birth position, $V(w_1)$], while all other
particles in the first generation, $\{x\dvtx |x|=1, x\not= w_1\}$ generate
independent copies of $\Theta$ (with respect to their birth positions).
The process goes on. The new system is denoted by $\mathcal{B}$.

\begin{fact}[(Lyons~\cite{lyons})]
\label{fLyons}
Assume (\ref{cond-hab}). The branching random walk
under $\Q$,
has the distribution of $\mathcal{B}$.
For any $|x|=n$, we have
%
\begin{equation}
\Q(w_n = x | \mathscr{F}_n) = {\ee^{-V(x)}\over W_n}.
\label{omega}
\end{equation}
The spine process $(V(w_n))_{n\ge0}$ under $\Q$
has the distribution of $(S_n)_{n\ge0}$ introduced
in Section \ref{snew-pair}.
\end{fact}

We mention that the analogue of Fact \ref{fLyons} for the branching
Brownian motion was known to Chauvin and Rouault~\cite{chauvin-rouault}.

Fact \ref{fLyons} is useful in the proof of the following
probabilistic estimate.

\begin{lemma}
\label{lminimal-position1}
Assume (\ref{cond-hab}), (\ref{cond-2}) and (\ref{cond-X}).
Let $C>0$ be the constant in Lemma \ref{la5}.
There exists a constant $c_{14}>0$ such that for
all sufficiently large $n$,
\[
\P \bigl\{ \exists x\dvtx n\le|x| \le2n, \tfrac{1}{2} \log n \le V(x)
\le\tfrac{1}{2} \log n + C \bigr\} \ge c_{14}.
\]
\end{lemma}

\begin{pf*}{Proof of Lemma \ref{lminimal-position1}} The
proof of
the lemma borrows an idea from \cite{elie}; see (\ref{elie}) below. We
fix $n$ and let
\[
a_i = a_i(n):= \cases{\displaystyle 0, &\quad $\mbox{if $\displaystyle0\le i\le
{n\over2}$,}$
\vspace*{2pt}\cr
\displaystyle{1\over2}\log n, &\quad $\mbox{if $\displaystyle
{n\over2}<i\le2n$}$}
\]
and for $n< k\le2n$,
\[
b_i^{(k)} = b_i^{(k)}(n):= \cases{
i^{1/12}, &\quad $\mbox{if $\displaystyle0\le i\le{n\over2}$,}$
\vspace*{2pt}\cr
(k-i)^{1/12}, &\quad $\mbox{if $\displaystyle{n\over2}<i\le k$.}$}
\]

For any vertex $y$, let, as before, $y_i$ denote the ancestor
of $y$ at generation $i$ (for $0\le i\le|y|$, with $y_{|y|}:= y$), and
$\Omega(y)$ the set of brothers of $y$. We consider
\begin{eqnarray*}
Z^{(n)} &:=& \sum_{k=n+1}^{2n}
Z^{(n)}_k,
\\
Z^{(n)}_{k} &:=& \# (E_k\cap F_k),
\end{eqnarray*}
where
\begin{eqnarray*}
E_k &:=& \bigl\{ y\dvtx |y|=k, V(y_i) \ge a_i,
\forall0\le i\le k, V(y) \le\tfrac{1}{2}\log n + C \bigr\},
\\
F_k &:=& \biggl\{ y\dvtx |y|=k, \sum_{v\in\Omega(y_{i+1})}
\bigl[1+\bigl(V(v)-a_{i}\bigr)^+\bigr] \ee^{-(V(v)-a_{i})} \le
c_{15}\ee^{-b_i^{(k)}},\\
&&\hspace*{218pt} \forall0\le i\le k-1 \biggr\}.
\end{eqnarray*}

[So if $x\in E_k$, then ${1\over2}\log n \le V(x) \le
{1\over2}\log n + C$. The set $E_k$ here has nothing to do with the
event $E_n$ in (\ref{En}).] The constant $c_{15}$ in the definition of
$F_k$ is positive and will be set later on. We make use of the new
probability measure $\Q$ introduced in (\ref{Q}): for $n<k\le2n$,
\[
\E\bigl[Z^{(n)}_k\bigr] = \E_\Q \biggl[
{Z^{(n)}_k \over W_k} \biggr] = \E_\Q \biggl[ \sum
_{|x|=k} {\mathbf{1}_{ \{ x\in E_k\cap F_k \} } \over
W_k} \biggr],
\]
which, by (\ref{omega}), is
\[
=\E_\Q\biggl[\sum_{|x|=k} \mathbf{1}_{
\{ x\in E_k\cap F_k \} } \ee^{V(x)} \mathbf{1}_{ \{ w_k =x \} } \biggr]
=\E_\Q\bigl[ \ee^{V(w_k)} \mathbf{1}_{ \{ w_k \in E_k\cap F_k \} } \bigr].
\]
Thus,
%
\begin{equation}
\label{eqznk} \E\bigl[Z^{(n)}_k\bigr] \ge n^{1/2}
\Q ( w_k \in E_k\cap F_k ).
\end{equation}

We need to estimate $\Q( w_k \in E_k\cap F_k )$. By Fact
\ref
{fLyons}, the process $(V(w_n))_{n\ge0}$ has the law of $(S_n)_{n\ge
0}$. Therefore, for $k\in(n, 2n] \cap\z$,
%
\begin{eqnarray}
\label{eqQEk} \Q ( w_k \in E_k ) &= &\P \biggl\{
S_i\ge a_i, \forall0\le i\le k, S_k \le
{1\over2}\log n + C \biggr\}
\nonumber
\\[-8pt]
\\[-8pt]
\nonumber
 &\in& \biggl[{c_{16}\over n^{3/2}},
{c_{17} \over n^{3/2}} \biggr],
\end{eqnarray}
by Lemmas \ref{la3} and \ref{la5}. We now use Lemma C.1 of
\cite{elie}, stating that for any $\varepsilon>0$, it is possible to
choose the constant $c_{15}$ (appearing in the definition of $F_k$)
sufficiently large such that for all large $n$,
%
\begin{equation}
\max_{k: n<k\le2n} \Q ( w_k \in E_k,
w_k \notin F_k ) \le {\varepsilon\over n^{3/2}}.
\label{elie}
\end{equation}

(The uniformity in $k\in(n, 2n] \cap\z$ is not stated in
\cite{elie}, but the same proof holds.) In particular, choosing
$\varepsilon:= {c_{16} \over2}$ [$c_{16}$ being in (\ref{eqQEk})]
leads to the existence of $c_{15}$ such that for all large $n$,
\[
\Q ( w_k \in E_k, w_k \in F_k )
\ge {c_{16} \over2 n^{3/2}}.
\]

 It follows from (\ref{eqznk}) that for all sufficiently
large $n$,
%
\begin{equation}
\E\bigl[Z^{(n)}\bigr] \ge \sum_{k=n+1}^{2n}
n^{1/2}{c_{16}\over2n^{3/2}} \ge c_{18}. \label{EZn}
\end{equation}

We now estimate the second moment of $Z^{(n)}$. By definition,
\[
\E \bigl[ \bigl(Z^{(n)}\bigr)^2 \bigr] = \sum
_{k=n+1}^{2n} \sum_{\ell=n+1}^{2n}
\E \bigl[ Z^{(n)}_k Z^{(n)}_{\ell} \bigr]
\le 2\sum_{k=n+1}^{2n} \sum
_{\ell=n+1}^k \E \bigl[ Z^{(n)}_k
Z^{(n)}_{\ell} \bigr].
\]

 Using again the probability $\Q$, we have for $n<\ell\le
k\le2n$,
\begin{eqnarray*}
\E \bigl[ Z^{(n)}_k Z^{(n)}_{\ell} \bigr]
&=& \E_\Q \biggl[ Z^{(n)}_{\ell}
{Z^{(n)}_k\over W_k} \biggr] = \E_\Q \biggl[
Z^{(n)}_{\ell} \sum_{|x|=k}
{\mathbf{1}_{ \{ x\in E_k\cap F_k
\} } \over W_k} \biggr] \\
&=&\E_\Q\bigl[ Z^{(n)}_{\ell}
\ee^{V(w_k)} \mathbf{1}_{ \{ w_k \in E_k\cap F_k
\} } \bigr]
\end{eqnarray*}
by (\ref{omega}), and thus is bounded by $\ee^C n^{1/2} \E
_\Q
[ Z^{(n)}_{\ell} \mathbf{1}_{ \{ w_k \in E_k\cap F_k \} } ]$. Therefore,
\[
\E \bigl[ \bigl(Z^{(n)}\bigr)^2 \bigr] \le 2
\ee^C n^{1/2} \sum_{k=n+1}^{2n}
\sum_{\ell=n+1}^k \E_\Q \bigl[
Z^{(n)}_{\ell} \mathbf{1}_{ \{ w_k \in E_k\cap F_k \} } \bigr].
\]

 We now estimate $\E_\Q[ Z^{(n)}_{\ell} \mathbf{1}_{ \{ w_k
\in
E_k\cap F_k \} }]$ on the right-hand side. It will be more convenient
to work with $ Y^{(n)}_\ell:= \sum_{|x|=\ell} \mathbf{1}_{\{x\in E_\ell
\}
}$ which is greater than $Z^{(n)}_\ell$. Decomposing the sum $Y_{\ell
}^{(n)}$ (for $n<\ell\le2n$) along the spine yields that
\[
Y^{(n)}_{\ell} = \mathbf{1}_{\{ w_\ell\in E_\ell\}} + \sum
_{i=1}^{\ell} \sum_{y\in
\Omega
(w_i)}
Y^{(n)}_\ell(y),
\]
where $\Omega(w_i)$ is, as before, the set of the brothers of
$w_i$, and $Y^{(n)}_\ell(y):= \# \{ x\dvtx |x| = \ell, x\ge y, x\in
E_\ell
\}$ the number of descendants $x$ of $y$ at generation $\ell$ such that
$x\in E_\ell$. By Fact \ref{fLyons}, the branching random walk
emanating from $y\in\Omega(w_i)$ has the same law under $\Q$ and under
$\P$. Therefore, conditioning on $\mathscr{G}_{\infty}:=\sigma\{
V(w_j), w_j, \Omega(w_j), (V(y))_{y\in\Omega(w_j)}, j\ge0\}$, we
have, for $y\in\Omega(w_i)$,
\[
\E_{\Q} \bigl[ Y^{(n)}_\ell| \mathscr{G}_{\infty}
\bigr] = \varphi_{i,\ell}\bigl(V(y)\bigr),
\]
where, for $r\in\r$,
\[
\varphi_{i,\ell}(r):= \E \biggl[ \sum_{|x|=\ell-i}
\mathbf{1}_{\{ r+V(x_j) \ge a_{j+i}, \forall
0\le j\le\ell-i, r+ V(x) \le{(1/2)} \log n + C\}} \biggr].
\]

Consequently,
\begin{eqnarray*}
\E \bigl[ \bigl(Z^{(n)}\bigr)^2 \bigr] &\le& 2
\ee^C n^{1/2} \sum_{k=n+1}^{2n}
\sum_{\ell=n+1}^k \Q \{ w_k \in
E_k\cap F_k, w_\ell\in E_\ell \}
\\
&&{} + 2\ee^C n^{1/2} \sum_{k=n+1}^{2n}
\sum_{\ell=n+1}^k \sum
_{i=1}^\ell \E_{\Q} \biggl[ {
\mathbf{1}}_{ \{ w_k \in E_k\cap F_k \} } \sum_{y\in\Omega(w_i)}
\varphi_{i,\ell}\bigl(V(y)\bigr) \biggr].
\end{eqnarray*}

 In the first double sum on the right-hand side, if $\ell=k$,
we simply argue that $\Q\{ w_k \in E_k\cap F_k, w_\ell\in E_\ell\}
\le\Q\{ w_k \in E_k\} \le{c_{17}\over n^{3/2}}$ [by (\ref
{eqQEk})], so that $\sum_{k=n+1}^{2n} \Q\{ w_k \in E_k\cap F_k, w_k
\in E_k\} \le\sum_{k=n+1}^{2n} {c_{17}\over n^{3/2}} = {c_{17}\over
n^{1/2}}$. This leads to
\begin{eqnarray*}
\E \bigl[ \bigl(Z^{(n)}\bigr)^2 \bigr] &\le& 2
\ee^C c_{17} + 2\ee^C n^{1/2} \sum
_{k=n+2}^{2n} \sum
_{\ell=n+1}^{k-1} \Q \{ w_k \in E_k
\cap F_k, w_\ell\in E_\ell \}
\\
&& {}+ 2\ee^C n^{1/2} \sum_{k=n+1}^{2n}
\sum_{\ell=n+1}^k \sum
_{i=1}^\ell \E_{\Q} \biggl[ {
\mathbf{1}}_{ \{ w_k \in E_k\cap F_k \} } \sum_{y\in\Omega(w_i)}
\varphi_{i,\ell}\bigl(V(y)\bigr) \biggr].
\end{eqnarray*}

 Recall from (\ref{EZn}) that $\E[ Z^{(n)} ] \ge c_{18}$.
Since $\P(Z^{(n)}>0)\ge{ \{ \E[Z^{(n)}] \}^2\over\E[(Z^{(n)})^2]}$,
the proof of Lemma \ref{lminimal-position1} is reduced to showing the
following estimates: for some constants $c_{19}>0$ and $c_{20}>0$ and
all sufficiently large $n$,
%
\begin{eqnarray}
\sum_{k=n+2}^{2n} \sum
_{\ell=n+1}^{k-1} \Q \{ w_k \in
E_k, w_\ell\in E_\ell \} &\le&
{c_{19}\over n^{1/2}}, \label{paley2}
\\
\sum_{k=n+1}^{2n} \sum
_{\ell=n+1}^k \sum_{i=1}^\ell
\E_{\Q} \biggl[ \mathbf{1}_{ \{ w_k \in E_k\cap F_k \} } \sum
_{y\in\Omega(w_i)} \varphi_{i,\ell}\bigl(V(y)\bigr) \biggr] &\le&
{c_{20}\over n^{1/2}}. \label{paley3}
\end{eqnarray}

Let us first prove (\ref{paley2}). By Fact \ref{fLyons}, for $n<
\ell
<k\le2n$,
\begin{eqnarray*}
&&\Q\{ w_k \in E_k, w_\ell\in E_\ell\}
\\
&&\qquad= \P \bigl\{ S_i\ge a_i, \forall0\le i\le k,
S_\ell\le\tfrac{1}{2}\log n + C, S_k \le
\tfrac{1}{2}\log n + C \bigr\}
\\
&&\qquad= \E \bigl\{ \mathbf{1}_{\{ S_i\ge a_i, \forall0\le i\le\ell,
S_\ell\le{1\over2}\log n + C\} } p_{k, \ell}(S_\ell) \bigr
\},
\end{eqnarray*}
where\footnote{Since $\ell> n$, we have, by definition, $a_i
= {1\over2}\log n$ for $i\ge\ell$.} $p_{k, \ell}(r):= \P\{
r+S_j\ge
{1\over2}\log n, \forall1\le j\le k-\ell, r+S_{k-\ell} \le{1\over
2}\log n + C\}$ (for $r\ge{1\over2}\log n$). Applying Lemma \ref
{la1} to $a:= r- {1\over2}\log n$ and $b:= 0$, we obtain, for $r\ge
{1\over2}\log n$,
\[
p_{k, \ell}(r) \le c_{21} {r- {(1/2)}\log n +1\over(k-\ell)^{3/2}},
\]
which leads to
\begin{eqnarray*}
&&\Q\{ w_k \in E_k, w_\ell\in E_\ell\}\\
&&\qquad\le {c_{21}\over(k-\ell)^{3/2}} \E \biggl\{ \mathbf{1}_{\{ S_i\ge a_i, \forall0\le i\le\ell,
S_\ell\le({1}/{2})\log n + C\} }
\biggl(S_\ell- {1\over2}\log n +1\biggr) \biggr\}
\\
&&\qquad\le {(C+1) c_{21}\over(k-\ell)^{3/2}} \P \biggl\{ S_i\ge a_i,
\forall0\le i\le\ell, S_\ell\le{1\over2}\log n + C \biggr
\}
\\
&&\qquad\le {(C+1) c_{21}\over(k-\ell)^{3/2}} {c_{22}\over n^{3/2}},
\end{eqnarray*}
the last inequality following from Lemma \ref{la3}. This
readily yields (\ref{paley2}).

It remains to check (\ref{paley3}). By (\ref{many-to-one}),
%
\begin{eqnarray}\label{phi1}
\qquad&&\varphi_{i,\ell}(r) \nonumber\\
&&\qquad= \E \bigl[ \ee^{S_{\ell-i}} {
\mathbf{1}}_{\{ r+ S_j \ge a_{j+i},
\forall0\le j\le\ell-i,
r+ S_{\ell-i} \le({1/2}) \log n + C\}} \bigr]
\\
&&\qquad\le n^{1/2}\ee^{C -r} \P \bigl[ r+ S_j \ge
a_{j+i}, \forall0\le j\le\ell-i, r+ S_{\ell-i} \le
\tfrac{1}{2} \log n + C \bigr].\nonumber
\end{eqnarray}

{F}rom here, we bound $\varphi_{i,\ell}(r)$ differently
depending on whether $i\le{n\over2}$ or $i > {n\over2}$.

\textit{First case}: $i\le{n\over2}$. By considering the $j=0$ term, we
get $\varphi_{i,\ell}(r) =0$ for $r<0$. For $r\ge0$, we have, by
(\ref
{phi1}) and Lemma \ref{la3},
%
\begin{eqnarray}\label{phi2}
\varphi_{i,\ell}(r) &\le& n^{1/2}\ee^{C -r}
c_{23} {r+1\over n^{3/2}}
\nonumber
\\[-8pt]
\\[-8pt]
\nonumber
 &=& {\ee^C c_{23}\over n}
\ee^{-r} (r+1),
\end{eqnarray}
so that writing $c_{24}:= \ee^C c_{23}$ and $\E_{\Q}
[k,i,\ell]:= \E_{\Q} [ \mathbf{1}_{ \{ w_k \in E_k \} } \sum_{y\in
\Omega
(w_i)} \varphi_{i,\ell}(V(y)) ]$ for brevity,
\begin{eqnarray*}
\E_{\Q} [k,i,\ell] &\le& {c_{24}\over n} \E_{\Q}
\biggl[ \mathbf{1}_{ \{ w_k \in E_k\cap F_k \} } \sum_{y\in\Omega(w_i)} {
\mathbf{1}}_{ \{ V(y) \ge0 \} } \ee^{-V(y)} \bigl(V(y)+1\bigr) \biggr]
\\
&\le& {c_{24}\over n} \E_{\Q} \biggl[ {
\mathbf{1}}_{ \{ w_k \in E_k\cap F_k \} } \sum_{y\in\Omega(w_i)} \ee^{-V(y)}
\bigl(V(y)^+ + 1\bigr) \biggr].
\end{eqnarray*}
By definition, we have $\sum_{y\in\Omega(w_i)}\ee^{-V(y)}
(V(y)^+ +1)\le c_{15}\ee^{-(i-1)^{1/12}}$ when $w_k \in F_k$. It
yields that
\[
\E_{\Q} [k,i,\ell] \le {c_{24}c_{15}\over n}\ee^{-(i-1)^{1/12}}
\Q(w_k \in E_k) \le {c_{24}c_{15}c_{17}\over n^{5/2}}
\ee^{-(i-1)^{1/12}}
\]
by (\ref{eqQEk}). As a consequence,
%
\begin{equation}
\sum_{k=n+1}^{2n} \sum
_{\ell=n+1}^k \sum_{1\le i\le{n/2}}
\E_{\Q} \biggl[ \mathbf{1}_{ \{ w_k \in E_k \} } \sum
_{y\in\Omega(w_i)} \varphi_{i,\ell}\bigl(V(y)\bigr) \biggr] \le
{c_{25}\over n^{1/2}}. \label{paley21}
\end{equation}

\textit{Second (and last) case}: ${n\over2}<i\le\ell$. This time, we
bound $\varphi_{i,\ell}(r)$ slightly differently. Let us go back to
(\ref{phi1}). Since $i> {n\over2}$, we have $a_{j+i}= {1\over2}\log
n$ for all $0\le j\le\ell-i$, thus $\varphi_{i,\ell}(r) =0$ for $r<
{1\over2}\log n$, whereas for $r\ge{1\over2}\log n$, we have, by
Lemma \ref{la1},
\[
\varphi_{i,\ell}(r) \le n^{1/2} \ee^{C-r}
{c_{26}\over(\ell-i+1)^{3/2}} \biggl(r- {1\over2}\log n +1\biggr).
\]

 This is the analogue of (\ref{phi2}); noting that the factor
${1\over n}$ becomes ${n^{1/2}\over(\ell-i+1)^{3/2}}$ now. {F}rom
here, we can proceed as in the first case: writing again $\E_{\Q}
[k,i,\ell]:= \E_{\Q} [ \mathbf{1}_{ \{ w_k \in E_k \} } \sum_{y\in
\Omega
(w_i)} \varphi_{i,\ell}(V(y)) ]$ for brevity, we have
\begin{eqnarray*}
\E_{\Q} [k,i,\ell] &\le& {c_{26}\ee^C n^{1/2}\over(\ell-i+1)^{3/2}}\\
&&{}\times  \E_{\Q}
\biggl[ \mathbf{1}_{ \{ w_k \in E_k\cap F_k \} } \sum_{y\in\Omega(w_i)}
\ee^{-V(y)} \biggl[\biggl(V(y)- {1\over2}\log n\biggr)^+ +1
\biggr] \biggr]
\\
&\le& {c_{26}\ee^C c_{15} n^{1/2}\over(\ell-i+1)^{3/2}} {\ee
^{-(k-i+1)^{1/12}} \over n^{1/2}}\Q(w_k\in
E_k)
\\
&\le& {c_{27} \over(\ell-i+1)^{3/2}n^{3/2}} \ee^{-(k-i+1)^{1/12}},
\end{eqnarray*}
where the last inequality comes from (\ref{eqQEk}). Consequently,
\[
\sum_{k=n+1}^{2n} \sum
_{\ell=n+1}^k \sum_{{n\over2}<i\le\ell}
\E_{\Q} \biggl[ \mathbf{1}_{ \{ w_k \in E_k \} } \sum
_{y\in\Omega(w_i)} \varphi_{i,\ell}\bigl(V(y)\bigr) \biggr] \le
{c_{28}\over n^{1/2}}.
\]

Together with (\ref{paley21}), this yields (\ref{paley3}),
and completes the proof of Lem\-ma~\ref{lminimal-position1}.
\end{pf*}

We have now all the ingredients for the proof of Theorem \ref
{tminimal-position}.

\begin{pf*}{Proof of Theorem \ref{tminimal-position}} Assume
(\ref
{cond-hab}), (\ref{cond-2}) and (\ref{cond-X}). Let $K>0$.

The system being super-critical, assumption (\ref{cond-hab})
ensures\break
$\P
\{ \min_{|x|=1} V(x) <0 \} >0$. Therefore, there exists an integer
$L=L(K)\ge1$ such that
\[
c_{29}:= \P \Bigl\{ \min_{|x|=L} V(x) \le-K \Bigr\}
>0.
\]

Let $n_k:= (L+2)^k$, $k\ge1$, so that $n_{k+1} \ge2n_k +L$, $\forall
k$. For any $k$, let
\[
T_k:= \inf \biggl\{ i\ge n_k\dvtx \min_{|x|=i}
V(x) \le{1\over2} \log n_k + C \biggr\},
\]
where $C>0$ is the constant in Lemma \ref
{lminimal-position1}. If $T_k <\infty$, let $x_k$ be such that $|x_k|
= T_k$ and that $V(x) \le{1\over2} \log n_k + C$. (If there are
several such $x_k$, any one of them will do the job, e.g., the
one with the smallest Harris--Ulam index.) Let
\[
G_k:= \{ T_k \le2n_k\} \cap \Bigl\{ \min
_{|y|=L} \bigl[V(x_ky)-V(x_k)\bigr] \le-K
\Bigr\},
\]
where $x_k y$ is the concatenation of the words $x_k$ and
$y$. For any pair of positive integers $j<\ell$,
%
\begin{equation}
\P \Biggl\{ \bigcup_{k=j}^\ell
G_k \Biggr\} = \P \Biggl\{ \bigcup_{k=j}^{\ell-1}
G_k \Biggr\} + \P \Biggl\{ \bigcap_{k=j}^{\ell-1}
G_k^c \cap G_\ell \Biggr\}.
\label{Borel-Cantelli}
\end{equation}

 On $\{ T_\ell<\infty\}$, we have
\[
\P\{ G_\ell| \mathscr{F}_{T_\ell} \} = \mathbf{1}_{ \{ T_\ell\le
2n_\ell\} }
\P \Bigl\{ \min_{|x|=L} V(x) \le-K \Bigr\} = c_{30} {
\mathbf{1}}_{
\{ T_\ell\le2n_\ell\} }.
\]

Since $\bigcap_{k=j}^{\ell-1} G_k^c$ is $\mathscr{F}_{T_\ell
}$-measurable, we obtain
\begin{eqnarray*}
\P \Biggl\{ \bigcap_{k=j}^{\ell-1}
G_k^c \cap G_\ell \Biggr\}& =& c_{30}
\P \Biggl\{ \bigcap_{k=j}^{\ell-1}
G_k^c \cap\{ T_\ell\le 2n_\ell \}
\Biggr\}\\
& \ge& c_{30} \P\{ T_\ell\le2n_\ell\} -
c_{30} \P \Biggl\{ \bigcup_{k=j}^{\ell
-1}
G_k \Biggr\}.
\end{eqnarray*}

Recall that $\P\{ T_\ell\le2n_\ell\} \ge c_{14}$ (Lemma
\ref
{lminimal-position1}; for large $\ell$, say $\ell\ge j_0$). Combining
this with (\ref{Borel-Cantelli}) yields that
\[
\P \Biggl\{ \bigcup_{k=j}^\ell
G_k \Biggr\} \ge(1-c_{30}) \P \Biggl\{ \bigcup
_{k=j}^{\ell-1} G_k \Biggr\} +
c_{14} c_{30}, \qquad j_0\le j<\ell.
\]

Iterating the inequality leads to
\begin{eqnarray*}
\P \Biggl\{ \bigcup_{k=j}^\ell
G_k \Biggr\} &\ge&(1-c_{30})^{\ell-j} \P\{
G_j\} + c_{14}c_{30} \sum
_{i=0}^{\ell-j-1} (1-c_{30})^i\\
& \ge&
c_{14}c_{30} \sum_{i=0}^{\ell-j-1}
(1-c_{30})^i.
\end{eqnarray*}

 This yields $\P\{ \bigcup_{k=j}^\infty G_k \} \ge c_{14}$,
$\forall j\ge j_0$. As a consequence,\break  $\P( \limsup_{k\to\infty} G_k)
\ge c_{14}$.

On the event $\limsup_{k\to\infty} G_k$, there are infinitely many
vertices $x$ such that $V(x) \le{1\over2} \log|x| + C -K$. Therefore,
\[
\P \biggl\{ \liminf_{n\to\infty} \biggl( \min_{|x|=n}
V(x) - {1\over
2}\log n \biggr) \le C -K \biggr\} \ge
c_{14}.
\]

 The constant $K>0$ being arbitrary, we obtain
\[
\P \biggl\{ \liminf_{n\to\infty} \biggl( \min_{|x|=n}
V(x) - {1\over
2}\log n \biggr) = -\infty \biggr\} \ge
c_{14}.
\]

Let $0<\varepsilon<1$. Let $J_1\ge1$ be an integer such that
$(1-c_{14})^{J_1} \le\varepsilon$. Under $\P^*$, the system survives
almost surely; so there exists a positive integer $J_2$ sufficiently
large such that $\P^* \{ \sum_{|x|=J_2} 1 \ge J_1 \}\ge1-\varepsilon$.
By applying what we have just proved to the sub-trees of the vertices
at generation $J_2$, we obtain
\[
\P^* \biggl\{ \liminf_{n\to\infty} \biggl( \min_{|x|=n}
V(x) - {1\over
2}\log n \biggr) = -\infty \biggr\} \ge1-
(1-c_{14})^{J_1} - \varepsilon \ge1-2\varepsilon.
\]
Sending $\varepsilon$ to $0$ completes the proof of Theorem
\ref{tminimal-position}.
\end{pf*}

Theorem \ref{tminimal-position} leads to the following result for the
lower limits of\break  $\min_{|x|=n} V(x)$, which was proved in \cite
{yzpolymer} under stronger assumptions (namely, $\E[(\sum_{|x|=1}1)^{1+\delta}] + \E[ \sum_{|x|=1} \ee^{-(1+\delta)V(x)}] +
\E[
\sum_{|x|=1} \ee^{ \delta V(x)}]<\infty$ for some $\delta>0$, and
(\ref
{cond-hab})). Recall that $\P^* (\bullet):= \P( \bullet| \mbox
{nonextinction})$.

\begin{theorem}
\label{tminimal-position-liminf=12}
Assume (\ref{cond-hab}), (\ref{cond-2}) and (\ref{cond-X}).
We have
\[
\liminf_{n\to\infty} {1\over\log n} \min
_{|x|=n} V(x) = {1\over2},\qquad\mbox{$\P^*$-a.s.}
\]

\end{theorem}

\begin{pf}
In view of Theorem \ref
{tminimal-position}, we
only need to check that\break  $\liminf_{n\to\infty} {1\over\log n} \min_{|x|=n} V(x) \ge{1\over2}$, $\P^*$-a.s.

Let $k>0$ and $a<{1\over2}$. By formula (\ref{many-to-one}) and in its
notation,
\begin{eqnarray*}
\E \biggl( \sum_{|x|=n} \mathbf{1}_{\{ \underline{V}(x) > -k\} } \mathbf{
1}_{\{ V(x) \le a \log n\} } \biggr) &=&\E \bigl( \ee^{S_n} \mathbf{1}_{\{ \underline{S}_n > -k\} }
\mathbf{ 1}_{\{
S_n \le a \log n\} } \bigr)
\\
&\le& n^a \P ( \underline{S}_n > -k, S_n
\le a \log n ),
\end{eqnarray*}
which, according to Lemma \ref{la1}, is bounded by a
constant multiple of $n^a {(\log n)^2\over n^{3/2}}$, and which is
summable in $n$ if $a< {1\over2}$. Therefore, as long as $a< {1\over
2}$, we have
\[
\sum_{n\ge1} \sum_{|x|=n} {
\mathbf{1}}_{\{ \underline{V}(x) > -k\} } \mathbf{ 1}_{\{ V(x) \le a \log n\} } < \infty,\qquad \mbox{$\P$-a.s.}
\]

By Biggins~\cite{biggins-lindley}, condition $\E( \sum_{|x|=1} \ee^{-V(x)} ) =1$
in (\ref{cond-hab}) implies that\break  $\inf_{|x|=n} V(x) \to\infty$, $\P^*$-a.s.;
thus $\inf_{|x|\ge0} V(x)
>-\infty$, $\P^*$-a.s. Consequently, $\liminf_{n\to\infty} {1\over
\log
n} \min_{|x|=n} V(x) \ge a$, $\P^*$-a.s., for any $a< {1\over
2}$.
\end{pf}

\section{Some questions}
\label{squestions}

Let $(V(x))$ be a branching random walk satisfying
(\ref
{cond-hab}), (\ref{cond-2}) and (\ref{cond-X}). Let, as before, $\P^*
(\bullet):= \P( \bullet| \mbox{nonextinction})$. Theorem \ref
{tminimal-position} tells us that $\liminf_{n\to\infty} [\min_{|x|=n}
V(x) - {1\over2}\log n ] = -\infty$, $\P^*$-a.s., but it does not give
us any quantitative information about how this ``$\liminf$'' expression
goes to $-\infty$. This leads to our first open question.

\begin{question}
\label{qrate-minimal-position}
Is there a deterministic sequence $(a_n)$ with
$\lim_{n\to\infty} a_n =\infty$ such that
\[
-\infty < \liminf_{n\to\infty} {1\over a_n} \biggl( \min
_{|x|=n} V(x) - {1\over2}\log n \biggr) <0,\qquad \mbox{$\P^*$-a.s.}?
\]
\end{question}

Our second question concerns the additive martingale $W_n$. In Theorem
\ref{tmain2}, we have proved that $\limsup_{n\to\infty} n^{1/2}W_n =
\infty$, $\P^*$-a.s., but the rate at which this ``$\limsup$'' goes to
infinity remains unknown.

\begin{question}
\label{qadditive-martingale}
Study the rate at which the upper limits of $n^{1/2} W_n$ go to infinity
$\P^*$-almost surely.
\end{question}

Questions \ref{qrate-minimal-position} and \ref{qadditive-martingale}
are obviously related via the inequality $W_n \ge \exp[-\min_{|x|=n}
V(x)]$. It is, however, not clear whether answering one of the
questions will necessarily lead to answering the other.

About the lower limits of $W_n$, we have a conjecture.

\begin{conjecture}
\label{conjmartingale-ratio}
We would have
\[
\liminf_{n\to\infty} n^{1/2} W_n = \biggl(
{2\over\pi\sigma^2} \biggr)^{ 1/2} D_\infty, \qquad\mbox{$
\P^*$-a.s.},
\]
where $\sigma^2:= \E[\sum_{|x|=1} V(x)^2 \ee^{-V(x)}]$.
\end{conjecture}

\section*{Acknowledgment}
We are grateful to Andreas Kyprianou for enlightenment on spinal decompositions.

%




\printaddresses


\begin{thebibliography}{35}

\bibitem{addario-berry-reed}
%
\begin{barticle}[mr]
\bauthor{\bsnm{Addario-Berry},~\bfnm{Louigi}\binits{L.}} \AND
\bauthor{\bsnm{Reed},~\bfnm{Bruce}\binits{B.}}
(\byear{2009}).
\btitle{Minima in branching random walks}.
\bjournal{Ann. Probab.}
\bvolume{37}
\bpages{1044--1079}.
\bid{doi={10.1214/08-AOP428}, issn={0091-1798}, mr={2537549}}
\bptok{imsref}%
\end{barticle}
%
\endbibitem

\bibitem{elie}
%
\begin{barticle}[auto:STB|2013/04/11|08:11:48]
\bauthor{\bsnm{A{\"{\i}}d{\'e}kon},~\bfnm{E.}\binits{E.}}
(\byear{2013}).
\btitle{Convergence in law of the minimum of a branching random
walk}.
\bjournal{Ann. Probab.}
\bvolume{41}
\bpages{1362--1426}.
\bid{mr={3098680}}
\bptok{imsref}%
\end{barticle}
%
\endbibitem

\bibitem{elie-bruno}
%
\begin{barticle}[mr]
\bauthor{\bsnm{A{\"{\i}}d{\'e}kon},~\bfnm{Elie}\binits{E.}} \AND
\bauthor{\bsnm{Jaffuel},~\bfnm{Bruno}\binits{B.}}
(\byear{2011}).
\btitle{Survival of branching random walks with absorption}.
\bjournal{Stochastic Process. Appl.}
\bvolume{121}
\bpages{1901--1937}.
\bid{doi={10.1016/j.spa.2011.04.006}, issn={0304-4149}, mr={2819234}}
\bptok{imsref}%
\end{barticle}
%
\endbibitem

\bibitem{ezsimple}
%
\begin{barticle}[mr]
\bauthor{\bsnm{A{\"{\i}}d{\'e}kon},~\bfnm{Elie}\binits{E.}} \AND
\bauthor{\bsnm{Shi},~\bfnm{Zhan}\binits{Z.}}
(\byear{2010}).
\btitle{Weak convergence for the minimal position in a branching random
walk: A
simple proof}.
\bjournal{Period. Math. Hungar.}
\bvolume{61}
\bpages{43--54}.
\bid{doi={10.1007/s10998-010-3043-x}, issn={0031-5303}, mr={2728431}}
\bptok{imsref}%
\end{barticle}
%
\endbibitem

\bibitem{berestycki-harris-kyprianou}
%
\begin{barticle}[mr]
\bauthor{\bsnm{Berestycki},~\bfnm{J.}\binits{J.}},
\bauthor{\bsnm{Harris},~\bfnm{S.~C.}\binits{S.~C.}} \AND
\bauthor{\bsnm{Kyprianou},~\bfnm{A.~E.}\binits{A.~E.}}
(\byear{2011}).
\btitle{Traveling waves and homogeneous fragmentation}.
\bjournal{Ann. Appl. Probab.}
\bvolume{21}
\bpages{1749--1794}.
\bid{doi={10.1214/10-AAP733}, issn={1050-5164}, mr={2884050}}
\bptok{imsref}%
\end{barticle}
%
\endbibitem

\bibitem{bertoin-rouault}
%
\begin{barticle}[mr]
\bauthor{\bsnm{Bertoin},~\bfnm{Jean}\binits{J.}} \AND
\bauthor{\bsnm{Rouault},~\bfnm{Alain}\binits{A.}}
(\byear{2005}).
\btitle{Discretization methods for homogeneous fragmentations}.
\bjournal{J. Lond. Math. Soc. (2)}
\bvolume{72}
\bpages{91--109}.
\bid{doi={10.1112/S0024610705006423}, issn={0024-6107}, mr={2145730}}
\bptok{imsref}%
\end{barticle}
%
\endbibitem

\bibitem{biggins-mart-cvg}
%
\begin{barticle}[mr]
\bauthor{\bsnm{Biggins},~\bfnm{J.~D.}\binits{J.~D.}}
(\byear{1977}).
\btitle{Martingale convergence in the branching random walk}.
\bjournal{J.~Appl. Probab.}
\bvolume{14}
\bpages{25--37}.
\bid{issn={0021-9002}, mr={0433619}}
\bptok{imsref}%
\end{barticle}
%
\endbibitem

\bibitem{biggins-lindley}
%
\begin{barticle}[mr]
\bauthor{\bsnm{Biggins},~\bfnm{J.~D.}\binits{J.~D.}}
(\byear{1998}).
\btitle{Lindley-type equations in the branching random walk}.
\bjournal{Stochastic Process. Appl.}
\bvolume{75}
\bpages{105--133}.
\bid{doi={10.1016/S0304-4149(98)00016-7}, issn={0304-4149}, mr={1629030}}
\bptok{imsref}%
\end{barticle}
%
\endbibitem

\bibitem{biggins10}
%
\begin{bincollection}[mr]
\bauthor{\bsnm{Biggins},~\bfnm{J.~D.}\binits{J.~D.}}
(\byear{2010}).
\btitle{Branching out}.
In \bbooktitle{Probability and Mathematical Genetics}
(\beditor{\binits{N.~H.}\bfnm{N.~H.} \bsnm{Bingham}}
\AND
\beditor{\binits{C.~M.}\bfnm{C.~M.} \bsnm{Goldie}}, eds.).
\bseries{London Mathematical Society Lecture Note Series}
\bvolume{378}
\bpages{113--134}.
\bpublisher{Cambridge Univ. Press}, \blocation{Cambridge}.
\bid{mr={2744237}}
\bptok{imsref}%
\end{bincollection}
%
\endbibitem

\bibitem{biggins-kyprianou96}
%
\begin{bincollection}[mr]
\bauthor{\bsnm{Biggins},~\bfnm{J.~D.}\binits{J.~D.}} \AND
\bauthor{\bsnm{Kyprianou},~\bfnm{A.~E.}\binits{A.~E.}}
(\byear{1996}).
\btitle{Branching random walk: {S}eneta--{H}eyde norming}.
In \bbooktitle{Trees ({V}ersailles, 1995)}
(\beditor{\binits{B.}\bfnm{B.} \bsnm{Chauvin} \betal{et al.}}, eds.).
\bseries{Progress in Probability}
\bvolume{40}
\bpages{31--49}.
\bpublisher{Birkh\"auser}, \blocation{Basel}.
\bid{mr={1439971}}
\bptok{imsref}%
\end{bincollection}
%
\endbibitem

\bibitem{biggins-kyprianou97}
%
\begin{barticle}[mr]
\bauthor{\bsnm{Biggins},~\bfnm{J.~D.}\binits{J.~D.}} \AND
\bauthor{\bsnm{Kyprianou},~\bfnm{A.~E.}\binits{A.~E.}}
(\byear{1997}).
\btitle{Seneta--{H}eyde norming in the branching random walk}.
\bjournal{Ann. Probab.}
\bvolume{25}
\bpages{337--360}.
\bid{doi={10.1214/aop/1024404291}, issn={0091-1798}, mr={1428512}}
\bptok{imsref}%
\end{barticle}
%
\endbibitem

\bibitem{biggins-kyprianou04}
%
\begin{barticle}[mr]
\bauthor{\bsnm{Biggins},~\bfnm{J.~D.}\binits{J.~D.}} \AND
\bauthor{\bsnm{Kyprianou},~\bfnm{A.~E.}\binits{A.~E.}}
(\byear{2004}).
\btitle{Measure change in multitype branching}.
\bjournal{Adv. in Appl. Probab.}
\bvolume{36}
\bpages{544--581}.
\bid{doi={10.1239/aap/1086957585}, issn={0001-8678}, mr={2058149}}
\bptok{imsref}%
\end{barticle}
%
\endbibitem

\bibitem{biggins-kyprianou05}
%
\begin{barticle}[mr]
\bauthor{\bsnm{Biggins},~\bfnm{J.~D.}\binits{J.~D.}} \AND
\bauthor{\bsnm{Kyprianou},~\bfnm{A.~E.}\binits{A.~E.}}
(\byear{2005}).
\btitle{Fixed points of the smoothing transform: The boundary case}.
\bjournal{Electron. J. Probab.}
\bvolume{10}
\bpages{609--631}.
\bid{doi={10.1214/EJP.v10-255}, issn={1083-6489}, mr={2147319}}
\bptok{imsref}%
\end{barticle}
%
\endbibitem

\bibitem{borovkov-foss}
%
\begin{barticle}[mr]
\bauthor{\bsnm{Borovkov},~\bfnm{A.~A.}\binits{A.~A.}} \AND
\bauthor{\bsnm{Foss},~\bfnm{S.~G.}\binits{S.~G.}}
(\byear{2000}).
\btitle{Estimates for the excess of a random walk over an arbitrary boundary
and their applications}.
\bjournal{Theory Probab. Appl.}
\bvolume{44}
\bpages{231--253}.
\bptok{imsref}%
\end{barticle}
%
\endbibitem

\bibitem{chauvin-rouault}
%
\begin{barticle}[mr]
\bauthor{\bsnm{Chauvin},~\bfnm{Brigitte}\binits{B.}} \AND
\bauthor{\bsnm{Rouault},~\bfnm{Alain}\binits{A.}}
(\byear{1988}).
\btitle{K{PP} equation and supercritical branching {B}rownian motion
in the
subcritical speed area. {A}pplication to spatial trees}.
\bjournal{Probab. Theory Related Fields}
\bvolume{80}
\bpages{299--314}.
\bid{doi={10.1007/BF00356108}, issn={0178-8051}, mr={0968823}}
\bptok{imsref}%
\end{barticle}
%
\endbibitem

\bibitem{feller}
%
\begin{bbook}[auto:STB|2013/04/11|08:11:48]
\bauthor{\bsnm{Feller},~\bfnm{W.}\binits{W.}}
(\byear{1971}).
\btitle{An Introduction to Probability Theory and Its Applications},
\bedition{2nd} ed.
\bpublisher{Wiley}, \blocation{New York}.
\bptok{imsref}%
\end{bbook}
%
\endbibitem

\bibitem{harris99}
%
\begin{barticle}[mr]
\bauthor{\bsnm{Harris},~\bfnm{Simon~C.}\binits{S.~C.}}
(\byear{1999}).
\btitle{Travelling-waves for the {FKPP} equation via probabilistic arguments}.
\bjournal{Proc. Roy. Soc. Edinburgh Sect. A}
\bvolume{129}
\bpages{503--517}.
\bid{doi={10.1017/S030821050002148X}, issn={0308-2105}, mr={1693633}}
\bptok{imsref}%
\end{barticle}
%
\endbibitem

\bibitem{heyde}
%
\begin{barticle}[mr]
\bauthor{\bsnm{Heyde},~\bfnm{C.~C.}\binits{C.~C.}}
(\byear{1970}).
\btitle{Extension of a result of {S}eneta for the super-critical
{G}alton--{W}atson process}.
\bjournal{Ann. Math. Statist.}
\bvolume{41}
\bpages{739--742}.
\bid{issn={0003-4851}, mr={0254929}}
\bptok{imsref}%
\end{barticle}
%
\endbibitem

\bibitem{yzpolymer}
%
\begin{barticle}[mr]
\bauthor{\bsnm{Hu},~\bfnm{Yueyun}\binits{Y.}} \AND
\bauthor{\bsnm{Shi},~\bfnm{Zhan}\binits{Z.}}
(\byear{2009}).
\btitle{Minimal position and critical martingale convergence in branching
random walks, and directed polymers on disordered trees}.
\bjournal{Ann. Probab.}
\bvolume{37}
\bpages{742--789}.
\bid{doi={10.1214/08-AOP419}, issn={0091-1798}, mr={2510023}}
\bptok{imsref}%
\end{barticle}
%
\endbibitem

\bibitem{jaffuel}
%
\begin{barticle}[auto:STB|2013/04/11|08:11:48]
\bauthor{\bsnm{Jaffuel},~\bfnm{B.}\binits{B.}}
(\byear{2012}).
\btitle{The critical barrier for the survival of the branching
random walk with absorption}.
\bjournal{Ann. Inst. H. Poincar\'e Probab. Statist.}
\bvolume{48}
\bpages{989--1009}.
\bid{mr={3052402}}
\bptok{imsref}%
\end{barticle}
%
\endbibitem

\bibitem{kahane-peyriere}
%
\begin{barticle}[mr]
\bauthor{\bsnm{Kahane},~\bfnm{J.~P.}\binits{J.~P.}} \AND
\bauthor{\bsnm{Peyri{\`e}re},~\bfnm{J.}\binits{J.}}
(\byear{1976}).
\btitle{Sur certaines martingales de {B}enoit {M}andelbrot}.
\bjournal{Adv. Math.}
\bvolume{22}
\bpages{131--145}.
\bid{issn={0001-8708}, mr={0431355}}
\bptok{imsref}%
\end{barticle}
%
\endbibitem

\bibitem{kozlov}
%
\begin{barticle}[mr]
\bauthor{\bsnm{Kozlov},~\bfnm{M.~V.}\binits{M.~V.}}
(\byear{1976}).
\btitle{The asymptotic behavior of the probability of non-extinction of
critical branching processes in a random environment}.
\bjournal{Theory Probab. Appl.}
\bvolume{21}
\bpages{791--804}.
\bptok{imsref}%
\end{barticle}
%
\endbibitem

\bibitem{kyprianou}
%
\begin{barticle}[mr]
\bauthor{\bsnm{Kyprianou},~\bfnm{A.~E.}\binits{A.~E.}}
(\byear{2004}).
\btitle{Travelling wave solutions to the {K}--{P}--{P} equation:
Alternatives to
{S}imon {H}arris' probabilistic analysis}.
\bjournal{Ann. Inst. Henri Poincar\'e Probab. Stat.}
\bvolume{40}
\bpages{53--72}.
\bid{doi={10.1016/S0246-0203(03)00055-4}, issn={0246-0203}, mr={2037473}}
\bptok{imsref}%
\end{barticle}
%
\endbibitem

\bibitem{lai}
%
\begin{barticle}[mr]
\bauthor{\bsnm{Lai},~\bfnm{Tze~Leung}\binits{T.~L.}}
(\byear{1976}).
\btitle{Asymptotic moments of random walks with applications to ladder
variables and renewal theory}.
\bjournal{Ann. Probab.}
\bvolume{4}
\bpages{51--66}.
\bid{mr={0391265}}
\bptok{imsref}%
\end{barticle}
%
\endbibitem

\bibitem{lalley-sellke}
%
\begin{barticle}[mr]
\bauthor{\bsnm{Lalley},~\bfnm{S.~P.}\binits{S.~P.}} \AND
\bauthor{\bsnm{Sellke},~\bfnm{T.}\binits{T.}}
(\byear{1987}).
\btitle{A conditional limit theorem for the frontier of a branching {B}rownian
motion}.
\bjournal{Ann. Probab.}
\bvolume{15}
\bpages{1052--1061}.
\bid{issn={0091-1798}, mr={0893913}}
\bptok{imsref}%
\end{barticle}
%
\endbibitem

\bibitem{liu98}
%
\begin{barticle}[mr]
\bauthor{\bsnm{Liu},~\bfnm{Quansheng}\binits{Q.}}
(\byear{1998}).
\btitle{Fixed points of a generalized smoothing transformation and applications
to the branching random walk}.
\bjournal{Adv. in Appl. Probab.}
\bvolume{30}
\bpages{85--112}.
\bid{doi={10.1239/aap/1035227993}, issn={0001-8678}, mr={1618888}}
\bptok{imsref}%
\end{barticle}
%
\endbibitem


\bibitem{lyons}
%
\begin{bincollection}[mr]
\bauthor{\bsnm{Lyons},~\bfnm{Russell}\binits{R.}}
(\byear{1997}).
\btitle{A simple path to {B}iggins' martingale convergence for
branching random
walk}.
In \bbooktitle{Classical and Modern Branching Processes
({M}inneapolis, {MN},
1994)}
(\beditor{\binits{K.~B.}\bfnm{K.~B.}~\bsnm{Athreya}}
\AND
\beditor{\binits{P.}\bfnm{P.}~\bsnm{Jagers}}, eds.).
\bseries{IMA Vol. Math. Appl.}
\bvolume{84}
\bpages{217--221}.
\bpublisher{Springer}, \blocation{New York}.
\bid{doi={10.1007/978-1-4612-1862-3_17}, mr={1601749}}
\bptok{imsref}%
\end{bincollection}
%
\endbibitem

\bibitem{lyons-peres}
%
\begin{bmisc}[auto:STB|2013/04/11|08:11:48]
\bauthor{\bsnm{Lyons},~\bfnm{R.}\binits{R.}}
\AND
\bauthor{\bsnm{Peres},~\bfnm{Yuval}\binits{Y.}}
(\byear{2010}).
\bhowpublished{\textit{Probability on Trees and Networks.}
Cambridge Univ. Press. Preprint. Available at \url
{http://mypage.iu.edu/\textasciitilde rdlyons/}}.
\bptok{imsref}%
\end{bmisc}
%
\endbibitem


\bibitem{lyons-pemantle-peres}
%
\begin{barticle}[mr]
\bauthor{\bsnm{Lyons},~\bfnm{Russell}\binits{R.}},
\bauthor{\bsnm{Pemantle},~\bfnm{Robin}\binits{R.}} \AND
\bauthor{\bsnm{Peres},~\bfnm{Yuval}\binits{Y.}}
(\byear{1995}).
\btitle{Conceptual proofs of {$L\log L$} criteria for mean behavior of
branching processes}.
\bjournal{Ann. Probab.}
\bvolume{23}
\bpages{1125--1138}.
\bid{issn={0091-1798}, mr={1349164}}
\bptok{imsref}%
\end{barticle}
%
\endbibitem

\bibitem{mckean75}
%
\begin{barticle}[mr]
\bauthor{\bsnm{McKean},~\bfnm{H.~P.}\binits{H.~P.}}
(\byear{1975}).
\btitle{Application of {B}rownian motion to the equation of
{K}olmogorov--{P}etrovskii--{P}iskunov}.
\bjournal{Comm. Pure Appl. Math.}
\bvolume{28}
\bpages{323--331}.
\bid{issn={0010-3640}, mr={0400428}}
\bptok{imsref}%
\end{barticle}
%
\endbibitem

\bibitem{mckean76}
%
\begin{barticle}[mr]
\bauthor{\bsnm{McKean},~\bfnm{H.~P.}\binits{H.~P.}}
(\byear{1976}).
\btitle{A correction to: ``{A}pplication of {B}rownian motion to the equation
of {K}olmogorov--{P}etrovski\u\i--{P}iskonov'' (\textit{{C}omm. {P}ure
{A}ppl. {M}ath.}
\textbf{28} (1975) 323--331)}.
\bjournal{Comm. Pure Appl. Math.}
\bvolume{29}
\bpages{553--554}.
\bid{issn={0010-3640}, mr={0423558}}
\bptok{imsref}%
\end{barticle}
%
\endbibitem

\bibitem{mogulskii73}
%
\begin{barticle}[mr]
\bauthor{\bsnm{Mogul'ski{\u\i}},~\bfnm{A.~A.}\binits{A.~A.}}
(\byear{1973}).
\btitle{Absolute estimates for moments of certain boundary functionals}.
\bjournal{Theory Probab. Appl.}
\bvolume{18}
\bpages{340--347}.
\bptok{imsref}%
\end{barticle}
%
\endbibitem

\bibitem{neveu86}
%
\begin{barticle}[mr]
\bauthor{\bsnm{Neveu},~\bfnm{J.}\binits{J.}}
(\byear{1986}).
\btitle{Arbres et processus de {G}alton--{W}atson}.
\bjournal{Ann. Inst. Henri Poincar\'e Probab. Stat.}
\bvolume{22}
\bpages{199--207}.
\bid{issn={0246-0203}, mr={0850756}}
\bptok{imsref}%
\end{barticle}
%
\endbibitem

\bibitem{seneta}
%
\begin{barticle}[mr]
\bauthor{\bsnm{Seneta},~\bfnm{E.}\binits{E.}}
(\byear{1968}).
\btitle{On recent theorems concerning the supercritical {G}alton--{W}atson
process.}
\bjournal{Ann. Math. Statist.}
\bvolume{39}
\bpages{2098--2102}.
\bid{issn={0003-4851}, mr={0234530}}
\bptok{imsref}%
\end{barticle}
%
\endbibitem

\end{thebibliography}
\end{document}